\documentclass[a4paper, fleqn]{amsart}
\usepackage[english]{babel}
\usepackage{lmodern}
\usepackage[T1]{fontenc}
\usepackage[latin1]{inputenc}
\usepackage{array}
\usepackage{amsfonts}
\usepackage{amssymb}
\usepackage{amsmath}
\usepackage{mathrsfs}
\usepackage{esint}
\usepackage[all]{xy}
\usepackage{enumitem}
\usepackage{a4wide}
\usepackage{concrete}
\usepackage{eulervm}
\usepackage{tikz}
\usepackage[bottom]{footmisc}

\numberwithin{equation}{section}
\newtheoremstyle{mytheoremstyle} 
    {4mm}                    
    {4mm}                    
    {\itshape}                   
    {6mm}                           
    {\scshape}                   
    {.}                          
    {0.5em}                       
    {}  

\theoremstyle{mytheoremstyle}

\newtheorem{df}{Definition}[section]
\let\olddf\df
\renewcommand{\df}{\olddf\normalfont}
\newtheorem{thm}[df]{Theorem}
\newtheorem{prop}[df]{Proposition}

\newtheorem{lem}[df]{Lemma}
\newtheorem{cor}[df]{Corollary}
\newtheorem{ex}[df]{Example}
\let\oldex\ex
\renewcommand{\ex}{\oldex\normalfont}
\newtheorem{rk}[df]{Remark}
\let\oldrk\rk
\renewcommand{\rk}{\oldrk\normalfont}
\newtheorem*{pr}{Proof}
\let\oldpr\pr
\renewcommand{\pr}{\oldpr\normalfont}

\newcommand{\C}{\mathbb{C}}
\newcommand{\R}{\mathbb{R}}
\newcommand{\Z}{\mathbb{Z}}
\newcommand{\N}{\mathbb{N}}
\newcommand{\F}{\mathcal{F}}
\newcommand{\A}{\mathcal{A}}

\renewcommand{\S}{\mathcal{S}}
\newcommand{\D}{\mathcal{D}}
\newcommand{\B}{\mathcal{B}}

\newcommand{\Tr}{\mathrm{Tr}}
\newcommand{\tr}{\mathrm{tr}}

\newcommand{\Ind}{\mathrm{Ind}}

\newcommand{\ch}{\mathrm{ch}}

\newcommand{\ad}{\mathrm{ad}}

\newcommand{\vol}{\mathrm{vol}}

\renewcommand{\dim}{\mathrm{dim}}

\renewcommand{\Re}{\mathrm{Re}}

\newcommand{\Res}{\mathrm{Res}}
\renewcommand{\vol}{\mathrm{vol}}
\newcommand{\ord}{\mathrm{ord}}
\newcommand{\dom}{\mathrm{dom}}

\renewcommand{\ch}{\mathrm{ch}}

\newcommand{\HP}{\mathrm{HP}}
\newcommand{\Pf}{\mathrm{Pf}}

\newcommand*{\barint}{\mathop{\ooalign{$\displaystyle{\int}$\cr$-$}}}
\newcommand*{\littlebarint}{\mathop{\ooalign{$\int$\cr$-$}}}
\renewcommand{\dom}{\mathrm{dom}}

\newcommand{\Op}{\mathrm{Op}}
\newcommand{\alg}{\mathrm{alg}}
\renewcommand{\\}{\vspace{2mm}}
\renewcommand{\i}{\textup{\textbf{i}}}
\newcommand{\CC}{\mathrm{CC}}
\newcommand{\Hom}{\mathrm{Hom}}
\newcommand{\Diff}{\mathrm{Diff}}
\renewcommand{\tilde}{\widetilde}
\renewcommand{\epsilon}{\varepsilon}

\title{Zeta Functions, Excision in Cyclic Cohomology and Index Problems}
\author{Rudy Rodsphon}
\address{Université de Lyon \\ CNRS UMR 5208 \\ Université Lyon 1 \\ Institut Camille Jordan \\  43, Bd du 11 novembre 1918, 69622 Villeurbanne Cedex, France}
\email{rodsphon@math.univ-lyon1.fr}

\begin{document}

\begin{abstract}
The aim of this paper is to show how zeta functions and excision in cyclic cohomology may be combined to obtain index theorems. In the first part, we obtain a local index formula for "abstract elliptic pseudodifferential operators" associated to spectral triples. This formula is notably well adapted when the zeta function has multiple poles. The second part is devoted to give a concrete realization of this formula by deriving an index theorem on the simple, but significant example of Heisenberg elliptic operators on a trivial foliation, which are in general non-elliptic but hypoelliptic. The last part contains a discussion on manifolds with conic singularity, more precisely about the regularity of spectral triples in this context.   \\

\setlength{\parindent}{0mm}

\end{abstract}

\maketitle 

\setlength{\parindent}{6mm}

\section*{Introduction}

Several years ago, Connes and Moscovici obtained in \cite{CM1995} a general index formula given in terms of residues of zeta functions, working with the so-called \emph{spectral triples}. A major advance was made since this formalism enlarges index theory to the more general context of the transverse geometry of foliations, where the interesting pseudodifferential operators are hypoelliptic without necessary being elliptic, while remaining Fredholm. Let us be a little more precise on this general formula. Connes and Moscovici constructed a Residue Cocycle on the algebra of the spectral triple, cohomologous to the Chern-Connes character in the $(B,b)$-complex of Connes. This cocycle has the feature of being "local", contrary to the representative of the Chern-Connes character obtained by changing the "Dirac operator" $D$ to the pseudodifferential operator $F = D \vert D \vert^{-1}$, which involves the operator trace, see \cite{ConIHES} or \cite{ConBook}. Here,  "local"  means that the cohomology class of the Residue Cocycle remains unchanged if the "Dirac operator" is perturbed by a smoothing operator. The interesting fact is that this happens because the Residue Cocycle is given by residues of zeta functions. Local index formulas are then deduced from a pairing between this cocycle and the K-theory of the algebra.    \\
 
In the spirit of the techniques developed by Connes and Moscovici, we give an abstract local index formula of a different flavour, which turns out to be useful to calculate the index of abstract elliptic pseudodifferential operators, in a sense to be defined. The formula is also given by a residue of a zeta function, but there is one important difference in that the cyclic cocycles concerned are defined not only on an "algebra of smooth functions", as in the Connes-Moscovici formula, but directly on the algebra of formal symbols of the pseudodifferential operators considered. We then illustrate on a simple but interesting example how such a formula may amount to topological index formulas, and in the end, discuss on the case of manifolds with conic singularity. Let us give an overview of the paper.  \\

Section 1 serves to recall some material about Higson's formalism developed in \cite{Hig2006}, concerning algebras of abstract differential operators and their relation with spectral triples, in particular regular ones. Following \cite{Uuy2009}, this allows to develop an abstract pseudodifferential calculus and a notion of ellipticity which covers many interesting examples. We shall focus on the example of Connes and Moscovici on foliations, involving the Heisenberg pseudodifferential calculus.  \\

The aim of Section 2 is to study the index theory in this context. More precisely, we construct a cyclic 1-cocycle on algebras of abstract pseudodifferential operators which generalizes the Radul cocycle defined for any closed manifold $M$, introduced by Radul in \cite{Rad1991}. The two important ingredients to construct this cocycle are, on the one hand, that the zeta function of a (classical) pseudodifferential operator on $M$ has a meromorphic extension to the complex plane, whose set of poles is at most simple and discrete. This allows the use of the Wodzicki-Guillemin residue. On the other hand, one uses the pseudodifferential extension and excision in periodic cyclic cohomology to push the trace on regularizing operators on $M$, viewed as a cyclic 0-cocycle, to a cyclic 1-cocycle on the algebra of formal symbols on $M$. The remarkable fact on using the Wodzicki-Guillemin residue is that it handles all the analytic issues, which will allow us to adopt an algebraic viewpoint in most of the paper. Excision in periodic cyclic cohomology then gives a local index formula for elliptic pseudodifferential operators, by compatibility with excision in K-theory.    

This construction is then extended to the abstract setting recalled in Section 1, and we obtain a cyclic 1-cocycle which generalizes the Radul cocycle in contexts where the zeta function exhibits multiple poles. 

\begin{thm} Let $\Psi(\Delta)$ be an algebra of abstract pseudodifferential operators on a Hilbert space $H$, and consider the pseudodifferential extension 
\[ 0 \to \Psi^{-\infty}(\Delta) \to \Psi(\Delta) \to \S = \Psi / \Psi^{-\infty} \to 0 \]
Suppose that the pole at zero of the zeta function is of order $p \geq 1$. Then, the cyclic 1-cocycle $\partial [\Tr] \in \HP^1(S)$, where $\Tr$ denotes the operator trace on $H$, is represented by the following functional, that we also call the Radul cocycle :
\begin{equation*} 
c(a_0,a_1) = \barint^1 a_0 \delta(a_1) - \dfrac{1}{2!}\barint^2 a_0 \delta^2(a_1) + \ldots + \dfrac{(-1)^{p-1}}{p!}\barint^p a_0 \delta^p(a_1) 
\end{equation*}
where $\delta(a) = [\log \Delta^{1/r}, a]$ and $\delta^k(a) = \delta^{k-1}(\delta (a))$ is defined by induction. The $r$ denotes the "order of $\Delta$"
\end{thm} 
The $\littlebarint^k$ are "higher Wodzicki-Guillemin residues" defined in Proposition \ref{higher WG trace}. \\

In Section 3, we show on an example how the results of the previous section may lead to index theorems, in the spirit of the Atiyah-Singer theorem. The example we work on is that of a trivial foliation $\R^p \times \R^q$, dealing with the Heisenberg pseudodifferential calculus. Even if this example is simple, it is also relevant for three reasons : Firstly, it allows to deal with hypoelliptic (non-elliptic) operators. Secondly, one can see how this leads to a purely algebraic approach of index theory ; analytic details are handled by the Wodzicki residue trace. Thirdly, the philosophy of the construction given is useful to understand how to adapt the techniques developed in \cite{Per2012} to treat for example the general case of foliations on closed manifolds (whose leaves are not necessarily compact). One interesting perspective is to obtain an index formula in the context of the transverse geometry of foliations, leading to a different approach as those of Connes and Moscovici in \cite{CM1998}. 

When dealing with the Radul cocycle, the main obstacle is that the formulas arising are, except in low dimensions, rather complicated. It is not obvious at all to obtain directly an index formula which depends only on the principal symbol. To cope with this difficulty, the general idea is to construct $(B,b)$-cocycles of higher degree which are cohomologous to the Radul cocycle in the $(B,b)$-bicomplex. These $(B,b)$-cocycles are shown to be more easily computable in the highest degree, for a reason that will be understood later. We give two ways of constructing these cocycles. In the first construction, we introduce homogeneous $(B,b)$-cocycles on regularizing operators, in many points similar to the cyclic cocycles associated to Fredholm modules given by Connes. The game still consists in pushing them to (inhomogeneous) $(B,b)$-cocycles on the algebra of Heisenberg (formal) symbols, using a zeta function regularisation of the trace and excision. The second construction involves Quillen's cochain theory in \cite{Qui1988}. The interest of using this formalism stands in the way we obtain the desired cocycles. Indeed, we do not have to go through the algebra of regularizing operators, so this method is completely algebraic. 

The context is a trivial foliation $\R^p \times \R^q$ of $\R^n$. Let $\S_H^0(\R^n)$ be the associated algebra of Heisenberg formal symbols of order $0$, and denote by 
\[ \sigma : \S_H^0(\R^n) \to C^\infty(S_H^* \R^n) \]
the principal symbol map. Here, $S_H^* \R^n$ denotes the "Heisenberg cosphere bundle", which is defined in Section \ref{CM pdo}.  Then, the main result of the section can be stated as follows : 

\begin{thm} The Radul cocycle is $(B,b)$-cohomologous to the homogeneous $(B,b)$-cocycle on $\S_H^0(\R^n)$ defined by 
\begin{equation*} 
\psi_{2n-1}(a_0, \ldots , a_{2n-1}) = \frac{(-1)^n}{(2\pi \i)^n} \int_{S_H^*\R^n} \sigma(a_0) d\sigma(a_1) \ldots d \sigma(a_{2n-1}) 
\end{equation*} 
\end{thm}

As an immediate corollary, we obtain the following index theorem.

\begin{thm} Let $P \in M_N(\Psi^0_H(\R^n))$ a Heisenberg elliptic pseudodifferential operator of formal symbol $u \in GL_N(\S^0_H(\R^n))$, and $[u] \in K_1(\S_H^0(\R^n))$ its (odd) K-theory class. Then, we have a formula for the Fredholm index of $P$ :
\[ \Ind(P) = \Tr(\Ind [u]) = - \frac{(n-1)!}{(2\pi \i)^n (2n-1)!} \int_{S_H^* \R^n} \tr(\sigma(u)^{-1} d\sigma(u) (d\sigma(u)^{-1} d\sigma(u))^{n-1})) \]
\end{thm}

Section 4 is a discussion on manifolds with conic singularity, and spectral triples associated. In this direction, note the work of Lescure in \cite{Lescure}, where spectral triples associated to conic manifolds are constructed. This construction has the notable feature that the zeta function associated has double order poles. The algebra considered in the spectral triple is the algebra of smooth functions vanishing to infinite order in a neighbourhood of the conic point, with a unit adjoined. Thus, many informations are lost in the differential calculus, e.g the abstract algebra of differential operators associated to the spectral triple cannot contain all the conic differential operators. Therefore, it is natural to ask if one can refine the choice of the algebra. Actually, we shall see that obtaining a regular spectral triple on such spaces inevitably leads us, in a certain manner, to erase the singularity. However, looking at this example gives a good picture of what happens when the regularity of the spectral triple is lost. The abstract Radul cocycle of Theorem 0.1, and thus the index formulas are no more local, because the terms killed by the residue in presence of regularity cannot be neglected in that case. We refer the reader to the concerned section for the different definitions and notations.

\begin{thm} Let $M$ be a conic manifold, i.e a manifold with boundary endowed with a conic metric, and let $r$ be a boundary defining function. Let $\Delta$ be the "conic laplacian" of Example \ref{conic laplacian}. Then, the Radul cocycle associated to the pseudodifferential extension  
\[ 0 \to r^{\infty} \Psi_b^{-\infty}(M) \to r^{-\Z} \Psi_b^{\Z}(M) \to r^{-\Z} \Psi_b^{\Z}(M)/r^{\infty} \Psi_b^{-\infty}(M) \to 0 \]
is given by the following \emph{non local} formula :
\begin{multline*}
c(a_0,a_1) = (\Tr_{\partial, \sigma} + \Tr_{\sigma})(a_0[\log \Delta, a_1]) - \frac{1}{2}\Tr_{\partial, \sigma}(a_0[\log \Delta,[\log \Delta, a_1]]) +  \\
+ \Tr_\partial \left( a_0 \sum_{k=1}^N a_1^{(k)} \Delta^{-k} \right) + \dfrac{1}{2\pi\i} \Tr\left(\int \lambda^{-z} a_0 (\lambda - \Delta)^{-1} a_1^{(N+1)} (\lambda - \Delta)^{-N-1}\right) \, d\lambda
\end{multline*} 
for $a_0, a_1 \in \Psi_b^{\Z}(M)/r^{\infty} \Psi_b^{-\infty}(M)$
\end{thm}

This approach yields another point of view on the eta invariant, the notable fact is that it is suitable also for pseudodifferential operators, and not only for Dirac operators. It might be an interesting problem to compare the formulas obtained with the usual eta invariant. \\

\textbf{Acknowledgements.} The author wishes to warmly thank Denis Perrot for sharing his insights, for his advices and constant support. He also thanks Thierry Fack for interesting discussions, and relevant remarks which helped to improve preliminary versions of this paper. The author is also grateful to Mathias Pétréolle for sharing some technical tips. 

\setlength{\parindent}{0mm}

\section{Abstract Differential Operators and Traces} \label{Abstract Differential Operators and a Trace}

In this part, we recall the Abstract Differential Operators formalism developed by Higson in \cite{Hig2006} to simplify the proof of the Connes-Moscovici local index formula \cite{CM1995}. This is actually another way of defining regular spectral triples. For details, the reader may refer to \cite{Hig2006} or \cite{Uuy2009}.

\subsection{Abstract Differential Operators} \label{Abstract Differential Operators}

Let $H$ be a (complex) Hilbert space and $\Delta$ a unbounded, positive and self-adjoint operator acting on it. To simplify matters, we suppose that $\Delta$ has a compact resolvent. \\

We denote by $H^\infty$ the intersection of all these Sobolev spaces. 
\begin{equation*}  
H^\infty = \bigcap_{k=0}^{\infty} \dom(\Delta^k)
\end{equation*}

\begin{df} \label{df ADO}
An algebra $\mathcal{D}(\Delta)$ of \emph{abstract differential operators} associated to $\Delta$ is an algebra of operators on $H^\infty$ fulfilling the following conditions
\begin{enumerate}[label={(\roman*)}]
\item The algebra $\D(\Delta)$ is filtered,
\begin{equation*} 
\D(\Delta) = \bigcup_{q=0}^{\infty} \D_q(\Delta)
\end{equation*}
\end{enumerate}
that is $\D_p(\Delta) \cdot \D_q(\Delta) \subset \D_{p+q}(\Delta)$. We shall say that an element $X \in \D_q(\Delta)$ is an \emph{abstract differential operator of order at most $q$}. The term \emph{differential order} will be often used for the order of such operators.  \\

\begin{enumerate}[label={(\roman*)}]
\item[(ii)] There is a $r>0$ ("the order of $\Delta$") such that for every $X \in \D_q(\Delta)$, $[\Delta,X] \in \D_{r+q-1}(\Delta)$. \\
\end{enumerate}

To state the last point, we define, for $s \in \R$, the $s$-\emph{Sobolev space} $H^s$ as the subspace $\dom(\Delta^{s/r})$ of $H$, which is a Hilbert space when endowed with the norm 
\[ \Vert v \Vert_s   = (\Vert v \Vert^2 + \Vert \Delta^{s/r}v \Vert^2)^{1/2} \]  

\begin{enumerate}[label={(\roman*)}]
\item[(iii)] \emph{Elliptic estimate.} If $X \in \D_{q}(\Delta)$, then, there is a constant $\varepsilon > 0$ such that 
\begin{equation*} 
\Vert v \Vert_q + \Vert v \Vert \geq \varepsilon \Vert Xv \Vert \, , \, \forall v \in H^\infty 
\end{equation*}
\end{enumerate}
\end{df}

Having Gärding's inequality in mind, the elliptic estimate exactly says that $\Delta^{1/r}$ should be thought as an "abstract elliptic operator" of order 1. It also says that any differential operator $X$ of order $q$ can be extended to a bounded operator form $H^{s+q}$ to $H^{s}$. This last property will be useful to define pseudodifferential calculus in this setting. \\

One example to keep in mind is the case in which $\Delta$ is a Laplace type operator on a closed Riemannian manifold $M$. Here, $r=2$ and $\D(\Delta)$ is simply the algebra of differential operators, the $H^s$ are the usual Sobolev space and we have an elliptic estimate. In fact, the definition above is an abstraction of this example, but it can be adapted to many more situations, for instance the case of foliations, on which we shall focus more in detail. 

\subsection{Correspondence with spectral triple}

Let $(A,H,D)$ a spectral triple (cf. \cite{CM1995} or \cite{Hig2006}). One may construct a algebra of abstract differential operators $\D = \D(A,D)$ inductively as follows : 
\begin{gather*} 
\D_0 = \text{algebra generated by } A \text{ and } [D,A] \\
\D_1 = [\Delta, \D_0] + \D_0[\Delta, \D_0] \\
\qquad \vdots \\
\D_k = \sum_{j=1}^{k-1} \D_j \cdot \D_{k-j} + [\Delta, \D_{k-1}] + \D_0[\Delta, \D_{k-1}]
\end{gather*}

Let $\delta$ be the unbounded derivation $\ad \vert D \vert = [\vert D \vert, \, . \,]$ on $\B(H)$. The spectral triple is $(A,H,D)$ is said \emph{regular} if $A, [D,A]$ are included in $\bigcap_{n=1}^{\infty} \dom \, \delta^{n}$. The following theorem of Higson makes the bridge between algebras of abstract differential operators and spectral triples.

\begin{thm} \emph{(Higson, \cite{Hig2006}).} Suppose that $A$ maps $H^\infty$ into itself. Then, the spectral triple $(A,H,D)$ is regular if and only if the elliptic estimate of Definition \ref{df ADO} holds. 
\end{thm}  

Regularity in spectral triples may be viewed an assumption allowing to control some asymptotic expansions of "pseudodifferential operators", as we shall see in the next paragraph from the perspective of the elliptic estimate.   

\subsection{Zeta Functions} \label{Zeta Functions}

Let $\D(\Delta)$ be an algebra of abstract differential operators. For $z \in \C$, one defines the \emph{complex powers} $\Delta^{-z}$ of $\Delta$ using functional calculus :
\begin{equation*}
\Delta^{-z} = \frac{1}{2 \pi i} \int \lambda^{-z} (\lambda - \Delta)^{-1} d\lambda
\end{equation*}
where the contour of integration is a vertical line pointing downwards separating $0$ and the (discrete) spectrum of $\Delta$. This converges in the operator norm when $\Re(z) > 0$, and using the semi-group property, all the complex powers can be defined after multiplying by $\Delta^k$, for $k \in \N$ large enough. Moreover, since $\Delta$ has compact resolvent, the complex powers of $\Delta$ are well defined operators on $H^\infty$. \\

We will suppose that there exists a $d\geq 0$ such that for every $X \in \D_q(\Delta)$, the operator $X \Delta^{-z}$ extends to a trace-class operator on $H$ for $z$ on the half-plane $\Re(z) > \frac{q+d}{r}$. The \emph{zeta function} of $X$ is 
\[ \zeta_X(z) = \Tr(X \Delta^{-z/r}) \]
The smallest $d$ verifying the above property is called the \emph{analytic dimension} of $\D(\Delta)$. In this case, the zeta function is holomorphic on the half-plane $\Re(z) > q+d$. We shall say that $\D(\Delta)$ has the \emph{analytic continuation property} if for every $X \in \D(\Delta)$, the associated zeta function extends to a meromorphic function of the whole complex plane. \\

There properties are set for all the section, unless if it is explicitly mentioned. \\

These notions come from properties of the zeta function on a closed Riemannian manifold $M$ : it is well-known that the algebra of differential operators on $M$ has analytic dimension $\dim \,M$ and the analytic continuation property. Its extension to a meromorphic function has at most simple poles at the integers smaller that $\dim \,M$. In the case where $M$ is foliated, the dimension of the leaves appears in the analytic dimension when working in the suitable context. Hence, the zeta function provide informations not only on the topology of $M$, but also on its the geometric structure,  illustrating the relevance of this abstraction.   

\subsection{Abstract Pseudodifferential Operators}  \label{Abstract Pseudodifferential Operators}

Let $\D(\Delta)$ an algebra of abstract differential operators of analytic dimension $d$. To define the notion of pseudodifferential operators, we need a more general notion of order, not necessary integral, which covers the one induced by the filtration of $\D(\Delta)$.

\begin{df}
An operator $T : H^\infty \to H^\infty$ is said to have \emph{pseudodifferential order} $m \in \R$ if for every $s \geq 0$, it extends to a bounded operator from $H^{m+s}$ to $H^s$. In addition, we require that operators of analytic order stricly less than $-d$ are trace-class operators.   
\end{df}

That this notion of order covers the differential order is due to the elliptic estimate, as already remarked in Section \ref{Abstract Differential Operators}. The space of such operators, denoted $\Op(\Delta)$, forms a $\R$-filtered algebra. There is also a notion of regularizing operators which are, as expected, the elements of the (two-sided) ideal of operators of all order. 

\begin{rk} Higson uses in \cite{Hig2006} the term "analytic order", but as the examples we deal with in the paper are about pseudodifferential operators, we prefer the term pseudodifferential order. 
\end{rk} 

\begin{ex} For every $\lambda \in \C$ not contained in the spectrum of $\Delta$, the resolvent $(\lambda - \Delta)^{-1}$ has analytic order $r$. Moreover, by spectral theory, its norm as an operator between Sobolev spaces is a $O(\vert \lambda \vert^{-1})$. 
\end{ex} 

The following notion is due to Uuye, cf. \cite{Uuy2009}. We just added an assumption on the zeta function which is necessary for what we do.  

\begin{df} \label{abstract pdo} An algebra of abstract pseudodifferential operators is a $\R$-filtered subalgebra $\Psi(\Delta)$ of $\Op(\Delta)$, also denoted $\Psi$ when the context is clear, satisfying 
\begin{align*}
& \Delta^{z/r} \Psi^m \subset \Psi^{\Re(z) + m},  \quad \Psi^m \Delta^{z/r} \subset \Psi^{\Re(z) + m}
\end{align*}
and which commutes, up to operators of lower order, with the complex powers of $\Delta^{1/r}$, that is , for all $m \in \R$, $z \in \C$
\begin{equation*} [\Delta^{z/r}, \Psi^m] \subset \Psi^{\Re(z) + m - 1} \end{equation*}

Moreover, we suppose that for every $P \in \Psi^m(\Delta)$, the zeta function 
\begin{equation*} 
\zeta_P(z) = \Tr(P \Delta^{-z/r}) 
\end{equation*}
is holomorphic on the half-plane $\Re(z) > m+d$, and extends to a meromorphic function of the whole complex plane. 
We shall denote by 
\begin{equation*} \Psi^{-\infty} = \bigcap_{m \in \R} \Psi^m \end{equation*}  
\end{df}


Of course, this is true for the algebra of (classical) pseudodifferential operators on a closed manifold. We shall recall later what happens in the example of Heisenberg pseudodifferential calculus on a foliation, as described by Connes and Moscovici in \cite{CM1995}.

We end this part with a notion of asymptotic expansion for abstract pseudodifferential operators. This can be seen as "convergence under the residue".

\begin{df} \label{asymptotic expansion} Let $T$ and $T_\alpha$ ($\alpha$ in a set $A$) be operators on $\Psi$. We shall write 
\begin{equation*} T \sim \sum_{\alpha \in A} T_\alpha \end{equation*} 
if there exists $c > 0$ and a finite subset $F \subset A$ such that for all finite set $F' \subset A$ containing $F$, the map
\[ z \longmapsto \Tr\left( (T - \sum_{\alpha \in F'} T_\alpha) \Delta^{z/r} \right)\]
is holomorphic in a half-plane $\Re(z) > -c$ (which contains $z=0$).
\end{df}

\begin{ex} Suppose that that for every $M > 0$, there exists a finite subset $F \subset A$ such that 
\begin{equation*} T - \sum_{\alpha \in F} T_\alpha \in \Psi^{-M} \end{equation*}
Then, $ T \sim \sum_{\alpha \in A} T_\alpha $
\end{ex}

In this context, asymptotic means that when taking values under the residue, such infinite sums, which have no reason to converge in the operator norm, are in fact finite sums. Thus, this will allow us to disregard analytic subtleties and to consider these sums only as formal expansions without wondering if they converge or not. In other words, this notion allows to adopt a completely algebraic viewpoint. To this effect, the following lemma is crucial. \\

\begin{lem} \label{CM trick} \emph{(Connes-Moscovici's trick, \cite{CM1995, Hig2006})} Let $Q \in \Psi(\Delta)$ be an abstract pseudodifferential operator. Then, for any $z \in \C$, we have 
\begin{equation}
[\Delta^{-z}, Q] \sim \sum_{k\geq 1} \binom{-z}{k} Q^{(k)}\Delta^{-z-k}
\end{equation}
where we denote $Q^{(k)} = \ad(\Delta)^k(Q)$, $\ad(\Delta) = [\Delta, \, . \,]$.
\end{lem}

\setlength{\parindent}{6mm}
\textsc{Two important facts.} 
\setlength{\parindent}{0mm}
Firstly, remark that the pseudodifferential order of terms in the sum are decreasing to $-\infty$, so that the difference between $[\Delta^{-z}, Q]$ and the sum becomes more and more regularizing as the number of terms grows.  

Secondly, and more importantly, this is essentially the elliptic estimate, or the regularity of the spectral triple, which implies this property. Then, if the sum in the lemma above is not asymptotic in the sense defined, the elliptic estimate cannot hold. In terms of spectral triple, this means it is not regular. 
 
\begin{pr} 
For $z \in \C$ of positive real part large enough, one proves, using Cauchy formulas and reasoning by induction, that the following identity holds (cf \cite{Hig2006}, Lemma 4.20) :
\begin{equation}  \label{CM trick bis}
\Delta^{-z} Q - Q \Delta^{-z} = \sum_{k = 1}^N \binom{-z}{k} Q^{(k)}\Delta^{-z-k} + \frac{1}{2 \pi i} \int \lambda^{-z}(\lambda - \Delta)^{-1} Q^{(N+1)}(\lambda - \Delta)^{-N-1} \, d\lambda
\end{equation}
By the elliptic estimate, the integral term in the right hand-side has pseudodifferential order $\ord Q +(N+1)r -N-1-(N+2)r = \ord(Q) - r - N -1$, which can therefore be made as small as we want by taking $N$ large. This proves the lemma in the case where $\Re(z) > 0$. The general case follows from the analytic continuation property. $\hfill{\square}$
\end{pr} 

\subsection{Higher traces on the algebra of abstract pseudodifferential operators} \label{Residue trace}

We give in this paragraph a simple generalization of the Wodzicki residue trace, when the zeta function of the algebra $\D(\Delta)$ has poles of arbitrary order. Actually, this was already noticed by Connes and Moscovici (see \cite{CM1995}). 

\begin{prop} \label{higher WG trace}
Let $\Psi(\Delta)$ an algebra of abstract pseudodifferential operators, following the context of the previous paragraphs. Suppose that the associated zeta function has a pole of order $p \geq 1$ in $0$. Then, the functional 
\begin{equation*} \barint^{p} P = \Res_{z=0} z^{p-1}\Tr(P \Delta^{-z/r})  \end{equation*}
defines a trace on $\Psi(\Delta)$.
\end{prop}

\begin{pr}
Let $P, Q \in \Psi(\Delta)$. Then, for $\Re(z) \gg 0$, we can use the trace property on commutators to write :
\begin{equation*} \Tr([P,Q] \Delta^{-z/r})  =  \Tr(P(Q - \Delta^{-z/r} Q \Delta^{z/r}) \Delta^{-z/r}) \end{equation*}
Hence, using the analytic continuation property, we have
\begin{equation*}\barint^{p} [P,Q]  =  \Res_{z=0} z^{p-1} \Tr(P(Q - \Delta^{-z/r} Q \Delta^{z/r}) \Delta^{-z/r}) \end{equation*}

By Lemma \ref{CM trick}, 
\begin{flalign*}  \Delta^{-z/r} Q - Q \Delta^{-z/r} \sim \sum_{k\geq 1} \binom{-z/r}{k} Q^{(k)}\Delta^{-k} \cdot \Delta^{-z/r}  \end{flalign*}
so that, 
\begin{flalign*}  \barint^p [P,Q] = \Res_{z=0} \sum_{k\geq 1} z^{p-1} \Tr \left( \binom{-z/r}{k} Q^{(k)}\Delta^{-k} \cdot \Delta^{-z/r} \right) \end{flalign*} 
The sum is finite : Indeed, the order of $Q^{(k)}\Delta^{-k}$ is $\ord(Q) - k$, so the terms in the sum above become holomorphic at $z=0$ when $k$ is large enough, and vanish when taking values under the residue. Finally, the finite sum remaining vanishes since the zeta function has at most a pole of order $p$ at $0$. $\hfill{\square}$
\end{pr}

If $0 \leq k < p$, then $\littlebarint^k$ is no more a trace in general, but one has an explicit relation expressing the commutators, cf. \cite{CM1995}. 

\subsection{The example of Connes and Moscovici} \label{CM pdo}

\subsubsection{Heisenberg pseudodifferential calculus on foliations} Let $M$ be a foliated manifold of dimension $n$, and let $\F$ be the integrable sub-bundle of the tangent bundle $TM$ of $M$ which defines the foliation. We denote the dimension of the leaves by $p$, and by $q = n-p$ their codimension.  \\

For the moment, we work in distinguished local charts. Let $(x_1, \ldots , x_n)$ a distinguished local coordinate system of $M$, i.e, the vector fields $\frac{\partial}{\partial x_1}, \ldots , \frac{\partial}{\partial x_p}$ (locally) span $\F$, so that $\frac{\partial}{\partial x_{p+1}}, \ldots , \frac{\partial}{\partial x_n}$ are transverse to the leaves of the foliation. Connes and Moscovici constructed in \cite{CM1995} an algebra of generalized differential operators using Heisenberg calculus, whose main idea is that : \\

\setlength{\parindent}{6mm}
\begin{itemize}
\item[$\bullet$] The vector fields $\frac{\partial}{\partial x_1}, \ldots , \frac{\partial}{\partial x_p}$ are of order 1. \\

\item[$\bullet$] The vector fields $\frac{\partial}{\partial x_{p+1}}, \ldots , \frac{\partial}{\partial x_n}$ are of order 2. \\
\end{itemize}
\setlength{\parindent}{0mm}

The Heisenberg pseudodifferential calculus consists in defining a class of smooth symbols $\sigma(x,\xi)$ on $\R^n_x \times \R^n_\xi$ which takes this notion of order into account. To this end, they set
\begin{align*}
& \vert \xi \vert' = (\xi_1^4 + \ldots + \xi_p^4 + \xi_{p+1}^2 + \ldots + \xi_{n}^2)^{1/4} \\ & \langle \alpha \rangle = \alpha_1 + \ldots + \alpha_p + 2\alpha_{p+1} + \ldots 2\alpha_n 
\end{align*}
for every $\xi \in \R^n$, $\alpha \in \N^n$.

\begin{df} A smooth function $\sigma(x,\xi) \in C^\infty(\R^n_x \times \R^n_\xi)$ is a Heisenberg symbol of order $m \in \R$ if $\sigma$ is $x$-compactly supported, and if for every multi-index $\alpha, \beta$, one has the following estimate
\begin{equation*} \vert \partial^\beta_x \partial^\alpha_\xi \sigma(x,\xi) \vert \leq (1+\vert \xi \vert')^{m - \langle \alpha \rangle} \end{equation*}
\end{df}

To such a symbol $\sigma$ of order $m$, one associates its left-quantization, which is the following operator 
\begin{equation*} 
P : C^\infty(\R^n) \to C^\infty(\R^n), \qquad Pf(x) = \frac{1}{(2\pi)^n} \int_{\R^n} e^{\i x \cdot \xi} \sigma(x,\xi) \hat{f}(\xi) d\xi 
\end{equation*}
We shall say that $P$ is a \emph{Heisenberg pseudodifferential operator of order $m$}, and denote the class of such operators by $\Psi_H^m(\R^n)$. The Heisenberg regularizing operators, whose class is denoted by $\Psi^{-\infty}(\R^n)$, are those of arbitrary order, namely 
\begin{equation*} \Psi^{-\infty}(\R^n) = \bigcap_{m \in \R} \Psi_H^m(\R^n) \end{equation*}
The reason why there is no $H$-subscript is that the Heisenberg regularizing operators are exactly the regularizing operators of the usual pseudodifferential calculus, i.e the operators with smooth Schwartz kernel.  \\

Actually we shall restrict to the smaller class of \emph{classical Heisenberg pseudodifferential operators}. For this, we first define the \emph{Heisenberg dilations}
\begin{equation*} 
\lambda \cdot (\xi_1, \ldots, \xi_p, \xi_{p+1}, \ldots, \xi_n) = (\lambda \xi_1, \ldots, \lambda \xi_p, \lambda^2 \xi_{p+1}, \ldots, \lambda^2 \xi_n) 
\end{equation*}
for any non-zero $\lambda \in \R$ and non-zero $\xi \in \R^n$. 

Then, a Heisenberg pseudodifferential operator $P \in \Psi^m_H(\R^n)$ of order $m$ is said \emph{classical} if its symbol $\sigma$ has an asymptotic expansion 
\begin{equation} \label{classical hpdo}
\sigma(x,\xi) \sim \sum_{j \geq 0} \sigma_{m-j}(x,\xi) 
\end{equation} 
where $\sigma_{m-j}(x,\xi) \in \S_H^{m-j}(\R^n)$ are \emph{Heisenberg homogeneous}, that is, for any non zero $\lambda \in \R$,
\begin{equation*} \sigma_{m-j}(x,\lambda \cdot \xi) = \lambda^{m-j} \sigma_{m-j}(x, \xi) \end{equation*}
The $\sim$ means that for every $M >0$, there exists an integer $N$ such that $\sigma - \sum_{j=0}^N \sigma_{m-j} \in \S_H^{-M}(\R^n)$.  \\
To avoid an overweight of notations, we shall keep the notation $\Psi_H$ to refer to classical elements. \\

Another important point is the behaviour of symbols towards composition of classical pseudodifferential operators. Of course, if $P,Q \in \Psi_H(\R^n)$ are Heisenberg pseudodifferential operators of symbols $\sigma_P$ and $\sigma_Q$, $PQ$ is also a Heisenberg pseudodifferential operator of order at most $\ord(P) + \ord(Q)$, and its the symbol $\sigma_{PQ}$ is given by the following asymptotic expansion called the \emph{star-product} of symbols, given by the formula
\begin{equation} \label{star product} 
\sigma_{PQ}(x, \xi) = \sigma_P \star \sigma_Q (x,\xi) \sim \sum_{\vert \alpha \vert \geq 0} \frac{(-\i)^{\vert \alpha \vert}}{\alpha!} \partial^\alpha_\xi \sigma_P(x, \xi) \partial^\alpha_x \sigma_Q(x, \xi) 
\end{equation} 
Note that the order of each symbol in the sum is decreasing while $\vert \alpha \vert$ is increasing. \\

We define the \emph{algebra of Heisenberg formal classical symbols} $\S_H(\R^n)$ as the quotient 
\begin{equation*}
\S_H(\R^n) = \Psi_H(\R^n)/\Psi^{-\infty}(\R^n) 
\end{equation*}
Its elements are formal sums given in (\ref{classical hpdo}), and the product is the star product (\ref{star product}). Note that the $\sim$ can be replaced by equalities when working at a formal level. \\ 
 
We now deal with ellipticity in this context. A Heisenberg pseudodifferential operator is said \emph{Heisenberg elliptic} if it is invertible in the unitalization  $\S_H(\R^n)^+$ of $\S_H(\R^n)$ . One can show that this is actually equivalent to say that its \emph{Heisenberg principal symbol}, e.g the symbol of higher degree in the expansion (\ref{classical hpdo}) is invertible on $\R^n_x \times \R^n_\xi \smallsetminus \{0\}$. An adaptation of arguments from classical elliptic regularity shows that the elliptic estimate holds in this case. A remarkable specificity of these operators is that they are hypoelliptic, but not elliptic in general. Nevertheless, they remain Fredholm operators between Sobolev spaces relative to this context. The interested reader should consult \cite{BeaGre} for details.    

\begin{ex} \label{subelliptic laplacian} The following operator, also called sub-elliptic sub-laplacian,   
\begin{equation*}
\Delta_H = \partial_{x_1}^4 +  \ldots + \partial_{x_p}^4 + \partial_{x_{p+1}}^2 + \ldots + \partial_{x_n}^2
\end{equation*} 
has Heisenberg principal symbol $\sigma(x,\xi) = \vert \xi \vert'^4$, and is therefore Heisenberg elliptic. However, its usual principal symbol, as an ordinary differential operator, is $(x,\xi) \mapsto \sum_{i=1}^p \xi_i^4 $, so $\Delta_H$ is clearly not elliptic. 
\end{ex}

Finally, Heisenberg pseudodifferential operators behaves well towards distinguished charts change. Therefore, Heisenberg pseudodifferential calculus can be defined globally on foliations by using a partition of unity. Then, for a foliated manifold $M$, we denote by $\Psi_H^m(M)$ the algebra of Heisenberg pseudodifferential operators on $M$. \\

It is not very difficult to verify the required assumptions of Definition \ref{abstract pdo}. However,  what concerns the zeta function is not obvious. 

\subsubsection{Residue Trace on Foliations} \label{Residue trace foliations} We now recall these results, proved by Connes and Moscovici in \cite{CM1995}.

\begin{thm} \emph{(Connes - Moscovici, \cite{CM1995})} Let $M$ be a foliated manifold of dimension $n$, $p$ be the dimensions of the leaves, and $P \in \Psi^m(M)$ be a Heisenberg pseudodifferential operator of order $m \in \R$. Let $\Delta$ the sub-elliptic sub-laplacian defined in Example \ref{subelliptic laplacian}, that we extend globally to $M$ by using a partition of unity. Then, the zeta function 
\begin{equation*} \zeta_P(z) = \Tr(P \Delta^{-z/4}) \end{equation*}
is holomorphic on the half-plane $\Re(z) > m+p+2q$, and extends to a meromorphic function of the whole complex plane, with at most simple poles in the set 
\begin{equation*} \{ m+p+2q, m+p+2q - 1, \, \ldots \} \end{equation*}
\end{thm}

\begin{rk} The analytic dimension of the algebra of Heisenberg differential operators is then $p+2q$. The $p$ is the dimension of the leaves, the $"2"$ is the degree of the vector fields transverse to them. 
\end{rk}

The meromorphic extension of the zeta function given by this theorem allows the construction of a Wodzicki-Guillemin trace on $\S_H(M) = \Psi_H(M)/\Psi^{-\infty(M)}$.

\begin{thm} \emph{(Connes - Moscovici, \cite{CM1995})} \label{CM residue} The Wodzicki residue functional 
\begin{equation*}  
\barint : \S_H(M) \longrightarrow \C , \quad P \longmapsto \Res_{z=0} \Tr(P \Delta^{-z/4}) 
\end{equation*}
is a trace. It is the unique trace on $\S_H(M)$, up to a multiplicative constant. Moreover, for $P \in \Psi_H(M)$, we have the following formula, only depending on the symbol $\sigma$ of $P$. 
\begin{equation}  \label{wodzicki residue symbol 0}
\barint P = \frac{1}{(2\pi)^n} \int_{S_H^*M} \iota_L \left(\sigma_{-(p+2q)}(x,\xi) \, \frac{ \omega^n}{n!}\right) 
\end{equation}
\end{thm}

Here, $S_H^*M$ is the Heisenberg cosphere bundle, that is, the sub-bundle 
\begin{equation*} 
S_H^*M = \{(x,\xi) \in T^*M \, ; \, \vert \xi \vert' = 1 \} 
\end{equation*} 
$L$ is the generator of the Heisenberg dilations, $\iota$ stands for the interior product and $\omega$ denotes the standard symplectic form on $T^*M$. 

\begin{rk} All these results still holds for Heisenberg pseudodifferential operators acting on sections of a vector bundle $E$ over $M$ : In this case, the symbol $\sigma_{-(p+2q)}(x,\xi)$ above is at each point $(x, \xi)$ an endomorphism acting on the fibre $E_x$, and (\ref{wodzicki residue symbol 0}) becomes :
\begin{equation*} 
\barint P = \frac{1}{(2\pi)^n} \int_{S_H^*M} \iota_L \left(\tr (\sigma_{-(p+2q)}(x,\xi)) \, \frac{ \omega^n}{n!}\right) 
\end{equation*}
where $\tr$ denotes the trace of endomorphisms. 
\end{rk}

\section{The Radul cocycle for abstract pseudodifferential operators}

\subsection{Abstract index theorems} \label{Abstract index}

We begin with another abstract setting. Let $A$ be an associative algebra over $\C$, possibly without unit, and $I$ an ideal in $A$. The extension
\begin{equation*} 0 \to I \to A \to A /I \to 0 \end{equation*}
gives rise to the following excision diagram, relating algebraic K-theory and periodic cyclic homology 
\begin{equation} \label{Abstract index diagram}
\xymatrix{
    K_1^{\alg}(A /I) \ar[r]^{\Ind} \ar[d]^{\ch_1} & K_0^{\alg}(I) \ar[d]^{\ch_0} \\
    \HP_1(A /I) \ar[r]^{\partial} & \HP_0(I)}
\end{equation}
The vertical arrows are respectively the odd and even Chern character. \\
 
We still denote $\partial :  \HP^0(I) \rightarrow \HP^1(A /I)$ the excision map in cohomology. As mentioned in \cite{Nis1997}, for $[\tau] \in \HP^0(I)$, $[u] \in K^1(A/I)$, one has the equality :
\begin{equation} \label{index pairing}
\langle [\tau], \ch_0 \Ind [u] \rangle = \langle \partial [\tau], \ch_1 [u] \rangle 
\end{equation}

One should have in mind that the left hand-side is an "analytic index", and think about the right hand-side as a "topological index".  \\

The construction of a boundary map $\partial$ in cohomology associated to an extension is standard. If $[\tau] \in \HP^0(I)$ is given by a hypertrace $\tau : I \to \C$, i.e a linear map satisfying the condition $\tau([A,I]) = 0$, then let us recall how to compute $\partial[\tau] \in  \HP^1(A /I)$. To begin, choose a lift $\tilde{\tau} : A \to \C$ of $\tau$, such that $\tilde{\tau}$ is linear (in general, this is not a trace), and a linear section $\sigma : A/I \to A$ such that $\sigma(1) = 1$, after adjoining a unit where we have to. Then, $\partial[\tau]$ is represented by the following cyclic cocycle :
\begin{equation*} 
c(a_0,a_1) = b\tilde{\tau}(\sigma(a_0), \sigma(a_1)) = \tilde{\tau}([\sigma(a_0), \sigma(a_1)])
\end{equation*} 
where $b$ is the Hochschild coboundary recalled in Section \ref{Bb complex}.

\subsection{The generalized Radul cocycle}

We can finally come to the main theorem of this section. Let $\D(\Delta)$ be an algebra of abstract differential operators and $\Psi$ be an algebra of abstract pseudodifferential operators. We consider the extension 
\begin{equation*} 0 \to \Psi^{-\infty} \to \Psi \to \S \to 0 \end{equation*}
where $\S$ is the quotient $\Psi / \Psi^{-\infty} $.
The operator trace on $\Psi^{-\infty}$ is well defined, and $\Tr([\Psi^{-\infty}, \Psi]) = 0$. 

\begin{thm} \label{local index formula} Suppose that the pole in zero of the zeta function is of order $p \geq 1$. Then, the cyclic 1-cocycle $\partial [\Tr] \in \HP^1(S)$ is represented by the following functional :
\begin{equation*} 
c(a_0,a_1) = \barint^1 a_0 \delta(a_1) - \dfrac{1}{2!}\barint^2 a_0 \delta^2(a_1) + \ldots + \dfrac{(-1)^{p-1}}{p!}\barint^p a_0 \delta^p(a_1) 
\end{equation*}
where $\delta(a) = [\log \Delta^{1/r}, a]$ and $\delta^k(a) = \delta^{k-1}(\delta(a))$ is defined by induction. We shall call this cocycle as the \emph{generalized Radul cocycle}. 
\end{thm} 

Here, the commutator $[\log \Delta^{1/r}, a]$ is defined as the non-convergent asymptotic expansion 
\begin{equation} \label{log commutator 1}
[\log \Delta^{1/r}, a] \sim \sum_{k\geq 1} \frac{(-1)^k}{k} a^{(k)} \Delta^{-k} 
\end{equation}
where $a^{(k)}$ has the same meaning as in Lemma \ref{CM trick}. 
This expansion arises by first using functional calculus :
\[ \log \Delta^{1/r} = \dfrac{1}{2 \pi \i} \int \log \lambda^{1/r} (\lambda - \Delta)^{-1} \, d\lambda \]
and then, reproducing the same calculations made in the proof of Lemma \ref{CM trick} to obtain the formula (cf. \cite{Hig2006} for details). In particular, note that $\log \Delta^{1/r} = \frac{1}{r} \log \Delta$.  \\

Another equivalent expansion possible, that we will also use, is the following
\begin{equation} \label{log commutator 2}
[\log \Delta^{1/r}, a] \sim \sum_{k\geq 1} \frac{(-1)^k}{k} a^{[k]} \Delta^{-k/r} 
\end{equation}
where $a^{[1]} = [\Delta^{1/r}, a]$, and $a^{[k+1]} = [\Delta^{1/r}, a^{[k]}]$. 
Before giving the proof of the result, let us give a heuristic explanation of how to get this formula. We first lift the trace on $\Psi^{-\infty}$ to a linear map $\tilde{\tau}$ on $\Psi$ using a zeta function regularization by "Partie Finie" :
\begin{equation*} \tilde{\tau}(P) = \Pf_{z=0} \Tr(P \Delta^{-z/r}) \end{equation*}
for any $P \in \Psi$. The "Partie Finie" $\Pf$ is defined as the constant term in the Laurent expansion of a meromorphic function. Let $Q \in \Psi$ be another pseudodifferential operator. Then, we have
\begin{equation*} \Pf_{z=0} \Tr([P,Q] \Delta^{-z/r}) = \Res_{z=0} \Tr \left(P \cdot \frac{Q - \Delta^{-z/r} Q \Delta^{z/r}}{z}\Delta^{-z/r} \right) \end{equation*}
by reasoning first for $z \in \C$ of sufficiently large real part to use the trace property, and then applying the analytic continuation property.
\setlength{\parindent}{0mm}
Then, informally we can think of the complex powers of $\Delta$ as 
\begin{equation*} \Delta^{z/r} = e^{\log \Delta \cdot z/r} = 1 + \dfrac{z}{r} \log \Delta + \ldots + \dfrac{1}{p!}\left(\dfrac{z}{r}\right)^p (\log \Delta)^{p} + O(z^{p+1}) \end{equation*}
which after some calculations, gives the expansion 
\begin{equation*}(Q - \Delta^{-z/r} Q \Delta^{z/r})\Delta^{-z/r} = z \delta(Q) - \dfrac{z^2}{2} \delta^2(Q) + \ldots + (-1)^{p-1}\dfrac{z^{p}}{p!}\delta^{p}(Q) + O(z^{p+1}) \end{equation*} 

\begin{pr}   Let $P,Q \in \Psi$ be two abstract pseudodifferential operators. The beginning of the proof is the same as the heuristic argument given above, so we start from the equality
\begin{align*} 
\Pf_{z=0} \Tr([P,Q] \Delta^{-z/r}) & = \Res_{z=0} \Tr \left(P \cdot \frac{Q - \Delta^{-z/r} Q \Delta^{z/r}}{z}\Delta^{-z/r} \right) \\
                                 & = \Res_{z=0} \Tr \left(P \cdot \dfrac{1}{z}\sum_{k\geq 1} \binom{-z/r}{k} Q^{(k)}\Delta^{-k} \cdot \Delta^{-z/r}\right)  
\end{align*} 
The second equality comes from Lemma \ref{CM trick}. \\

Then, let $X$ be an indeterminate. As power series over the complex numbers with indeterminate $X$, we remark that for any $z \in \C$, one has 
\begin{equation*} 
\dfrac{1}{z} \sum_{k\geq 1} \binom{-z/r}{k} X^k = \frac{1}{z}((1+X)^{-z/r} - 1) 
\end{equation*}

On the other hand, we have, for $q \in \N$,
\begin{equation*}
\ad(\log \Delta^{1/r})^{q}(Q) = \dfrac{1}{r^q}[\log \Delta , [... , [\log \Delta, Q] ] ] \sim \dfrac{1}{r^q}\sum_{k \geq q} \sum_{k_1 + \ldots + k_q = k} \frac{(-1)^k}{k_1 \ldots k_q} Q^{(k)} \Delta^{-k} 
\end{equation*}
Using once more the indeterminate $X$, one has
\begin{align*} 
\sum_{k \geq q} \sum_{k_1 + \ldots + k_q = k} \frac{(-1)^k}{k_1 \ldots k_q} X^k &= \left(\sum_{l \geq 1} \frac{(-1)^l X^l}{l} \right) \\
&= \log(1+X)^q
\end{align*}
thus obtaining 
\begin{equation*}  
\sum_{q \geq 1} \frac{(-1)^{q-1}}{q!}\dfrac{z^{q-1}}{r^q} \log(1+X)^q = \frac{1}{z}((1+X)^{-z/r} - 1) 
\end{equation*}
This proves that the coefficients of $Q^{(k)}\Delta^{-k}$ in the sums 
\begin{align*}
&\dfrac{1}{z}\sum_{k\geq 1} \binom{-z}{k} Q^{(k)}\Delta^{-k} , &
\sum_{q \geq 1} \frac{(-1)^{q-1}}{q!}  \dfrac{z^{q-1}}{r^q}\left(\sum_{k \geq q} \sum_{k_1 + \ldots + k_q = k} \frac{(-1)^k}{k_1 \ldots k_q} Q^{(k)} \Delta^{-k} \right) 
\end{align*}
are the same, hence the result follows. $\hfill{\square}$ 
\end{pr}

Applying the pairing (\ref{index pairing}), we have a local index theorem.

\begin{ex} \label{radul foliation}
Let $M$ be a closed foliated manifold with integrable sub-bundle $F \in TM$, $\Delta$ the sub-elliptic sub-laplacian of Example \ref{subelliptic laplacian} sand take $\Psi(\Delta) = \Psi_H(M)$ the algebra of (classical) Heisenberg pseudodifferential operators on $M$, $\Psi^{-\infty}(\Delta) = \Psi^{-\infty}(M)$ the ideal of regularizing operators. The quotient $\Psi/\Psi^{-\infty}$ is the algebra $\S_{H}(M)$ of full classical Heisenberg symbols. 
A trace on $\Psi^{-\infty}(M)$ is given by 
\begin{equation}  \label{usual trace}
\tau(K) = \Tr(K) = \int_M k(x,x) d\vol(x) 
\end{equation}
where $k$ is the Schwartz kernel of $K$. 
Then, using the residue defined in Theorem \ref{CM residue} and applying Theorem \ref{local index formula}, $\partial[\tau]$ is represented by the following cyclic 1-cocycle on $\S_H(M)$ :
\begin{equation} \label{radul cocycle}
c(a_0,a_1) = \barint a_0 [\log \vert \xi \vert', a_1] 
\end{equation}
With a slight abuse of notation, we denoted by $\log \vert \xi \vert'$ the symbol of $\Delta^{1/4}$. We emphasize that the product of symbols is the star-product defined in (\ref{star product}), but we omit the notation $\star$. \\

Remark that $\log \vert \xi \vert'$ is a log-polyhomogeneous (Heisenberg) symbol and is not classical. But using (\ref{log commutator 2}), it is clear that its commutator with any element of $\S_H(M)$ is. Note also that the cocycle is defined on the symbols rather that on the operators, but this does not matter since the residue kills the smoothing contributions. In particular, only a finite number of terms of the star-product are involved. This is exactly what we meant when we said that the Wodzicki residue handles analytic issues in the introduction.  \\

This cocycle was first introduced by Radul in \cite{Rad1991} in the context of closed manifold, without considering foliations, as a 2-cocycle over the Lie algebra of formal symbols on the manifold. The Radul cocycle also may appear from a Partie Finie regularization of the zeta function, so we keep the same name for the cocycle \ref{radul cocycle} obtained in this more general setting. \\

From this cocycle, we then get an index formula for Heisenberg elliptic pseudodifferential operators. Indeed, if $P$ is such an operator of formal symbol $u \in \S_H(M)$, and $Q$ a parametrix of $P$ in the Heisenberg calculus, of formal symbol $u' \in \S_H(M)$, then, the Fredholm index  of $P$ is given by 
\[ \Ind(P) = c(u,u') \]
As we can see, the Radul cocycle is given by a Wodzicki residue, and is hence local. However, it seems to be an unattainable task to get an index formula in terms of the principal symbol only since by (\ref{wodzicki residue symbol 0}), we have to find the symbol of order $-(p+2q)$ of $u [\log \vert \xi \vert', u']$. At first sight, many terms of the formal expansions of $u$ and $u'$, as well as many of their higher derivatives, seem to be involved. We shall see in next section a way to overcome this difficulty. 
\end{ex}

\section{A computation of the Radul cocycle} \label{computation}

At first sight, the latter index formula obtained is local in the sense that it is given as a residue formula, a little in the spirit of that of Connes and Moscovici. However, as already noted in Example \ref{radul foliation}, the formula obtained is rather involved. \\

This section is devoted to show how one may recover an interesting index formula from the Radul cocycle, working on the simplest foliation possible. For all this section, even if it is not explicitly mentioned, we consider $\R^n$ as a trivial foliation $\R^p \times \R^q$, where $0 \leq p < n$ and $q=n-p$, and consider the associated classical Heisenberg pseudodifferential operators $\Psi_H^0(\R^n)$ of order 0. \\

Our goal is to show that the Radul cocycle (\ref{radul cocycle}) on $\S_H^0(\R^n)$ is cohomologous in $\HP^1(\S_H(\R^n))$ to simple inhomogeneous $(B,b)$-cocycles of higher degree, making the computation of the index problem easier. We shall always use coordinates adapted to the foliation $\R^p \times \R^q$. \\

We shall give two ways of constructing these cocycles. Before beginning these constructions, we briefly recall how to define the $(B,b)$-bicomplex. 

\subsection{The $(B,b)$-bicomplex} \label{Bb complex} Let $A$ be an associative algebra over $\C$. For $k \geq 0$, denote by $\CC^{k}(A)$ the space of $(k+1)$-linear forms on the unitalization $A^+$ of $A$ such that $ \phi(a_0, \ldots, a_k) = 0 $ when $a_i = 1$ for some $i \geq 1$. Then, define the differentials 
\begin{equation*} B : \CC^{k+1}(A) \to \CC^{k}(A), \qquad b : \CC^{k}(A) \to \CC^{k+1}(A) \end{equation*}
by the formulas 
\begin{multline*}
B\phi(a_0, \ldots, a_k) = \sum_{i=0}^{k}(-1)^{ik} \phi(1, a_i, \ldots, a_k, a_0, \ldots, a_{i-1}) \\
b\phi(a_0, \ldots, a_{k+1}) = \sum_{i=0}^{k} (-1)^{i} \phi(a_0, \ldots, a_{i-1}, a_i a_{i+1}, a_{i+2}, \ldots, a_{k+1}) \\
+ (-1)^{k+1} \phi(a_{k+1} a_0, \ldots, a_k)
\end{multline*}
that is, $B^2 = b^2 = 0$. Moreover, $B$ and $b$ anticommute, which allows to define the $(B,b)$-bicomplex 
\begin{equation*}
\xymatrix{
     & \vdots  & \vdots  & \vdots & \\
    \ldots  \ar[r]^-{B} & \CC^2(A) \ar[r]^-{B} \ar[u]^-{b} & \CC^1(A) \ar[r]^-{B} \ar[u]^-{b} & \CC^0(A) \ar[u]^-{b} & \\
    \ldots  \ar[r]^-{B}  & \CC^1(A) \ar[r]^-{B} \ar[u]^-{b} & \CC^0(A) \ar[u]^-{b} & & \\
    \ldots  \ar[r]^-{B}  & \CC^0(A) \ar[u]^-{b} & & &}
\end{equation*}
Then, the periodic cyclic cohomology $\HP^{\bullet}(A)$ is the cohomology of the total complex. More precisely, it is the cohomology of the 2-periodic complex
\[ \xymatrix{\ldots \ar[r]^-{B+b} & \CC^\mathrm{even}(A) \ar[r]^-{B+b} & \CC^\mathrm{odd}(A) \ar[r]^-{B+b} & \CC^\mathrm{even}(A) \ar[r]^-{B+b} & \ldots } \]
where
\begin{gather*} 
\CC^\mathrm{even}(A) = \CC^{0}(A) \oplus \CC^{2}(A) \oplus \ldots \\
\CC^\mathrm{odd}(A) = \CC^{1}(A) \oplus \CC^{3}(A) \oplus \ldots
\end{gather*}
Hence, there are only an even and an odd periodic cyclic cohomology groups, respectively denoted $\HP^{0}(A)$ and $\HP^{1}(A)$.

\begin{rk} Sometimes, authors consider the total differential $B-b$ instead of $B+b$.  
\end{rk}

\subsection{General context} \label{general context} Recall from Section \ref{Residue trace} that the residue trace of a Heisenberg pseudodifferential operator $P \in \Psi_H(\R^n)$ of symbol $\sigma$ is given by 
\begin{equation} \label{wodzicki residue symbol}
\barint P = \frac{1}{(2\pi)^n} \int_{S_H^*\R^n} \iota_L \left(\sigma_{-(p+2q)}(x,\xi) \frac{\omega^n}{n!} \right)
\end{equation}
where $\sigma_{-(p+2q)}$ is the Heisenberg homogeneous term of order $-(p+2q)$ in the asymptotic expansion of $\sigma$, $\omega = \sum_i dx_i d\xi_i$ is the standard symplectic form on $T^*\R^n = \R^n_x \times \R^n_\xi$, and $L$ is the generator of the Heisenberg dilations, given by the formula 
\[ L = \sum_{i=1}^p \xi_i \partial_{\xi_i} + 2\sum_{i=p+1}^n \xi_i \partial_{\xi_i} \]  
Note that in this example, the sub-elliptic sub-laplacian has not a compact resolvent since we work on $\R^n$. However, the results in Section \ref{Residue trace foliations} on the Wodzicki residue still holds because we consider pseudodifferential operators which have compact support. \\

We first extend the trace on $\Psi^{-\infty}(\R^n)$ given in (\ref{usual trace}) to a graded trace on the graded algebra $\Psi^{-\infty}(\R^n) \otimes \Lambda^\bullet T^*\R^n$, using a Berezin integral :
\begin{equation*} 
\Tr(K \otimes \alpha) = \alpha_{[2n]}\Tr(K) 
\end{equation*}
where $K \in \Psi^{-\infty}(\R^n)$, and $\alpha_{[2n]}$ is the coefficient of the form $dx_1 \ldots dx_n d\xi_1 \ldots  d\xi_n$ in $\alpha$ (the wedges are dropped to simplify notations). Here, we emphasize once more that $T^*\R^n$ is seen as the vector space $\R^n_x \times \R^n_\xi$. Therefore $\Lambda^\bullet T^*\R^n$ stands for the exterior algebra of the vector space $T^*\R^n = \R^n_x \times \R^n_\xi$, and not for the vector bundle of exterior powers of the cotangent bundle, as usual.  \\

Moreover, the Wodzicki residue trace on $\Psi_H(\R^n)$ is given by a zeta function regularisation of this trace. Therefore, the latter procedure also extends the Wodzicki residue trace to a graded trace on the graded algebra $\Psi_H(\R^n) \otimes \Lambda^\bullet T^*\R^n$. The latter descends to a graded trace on $\S_H(\R^n)\otimes \Lambda^\bullet T^*\R^n$. The composition law of pseudodifferential operators, or the star-product of symbols for the latter, are extended to these algebras just by imposing that they commute to elements of the exterior algebra. \\

Remark also that the following commutation relations hold
\begin{align*}
&[x_i, \xi_j] = \i  \delta_{i,j}, \quad  [x_i,  x_j] = [\xi_i, \xi_j] = 0
\end{align*}
where we denote $\i = \sqrt{-1}$. In short, $\ad(x_i)$ and $\ad(\xi_i)$ are respectively the differentiation of symbols with respect to the variables $\xi_i$ and $x_i$.

Finally, let $F$ be the multiplier on $\S_H(\R^n)\otimes \Lambda^\bullet T^*\R^n$ defined by 
\begin{equation*} 
F = \sum_i (x_i d\xi_i + \xi_i dx_i) 
\end{equation*}
As the two following lemmas might indicate, this operator will play a role rather similar to operators usually denoted by $F$ when dealing with finitely summable Fredholm modules. The difference is that this $F$ here is not the main object of study, and acts more as an intermediate towards the main result. 

\begin{lem} $F^2$ is equal to $\i  \omega$, where $\omega$ is the standard symplectic form on $T^*\R^n$. In particular, $F^2$ commutes to every element in $\S_H(\R^n) \otimes \Lambda^\bullet T^*\R^n$. 
\end{lem}

\begin{lem} \label{F differential} For every symbol $a \in \S_H(\R^n)$, one has 
\begin{equation*} [F,a] = \i  da = \i  \sum_i \left( \frac{\partial a}{\partial x_i} dx_i + \frac{\partial a}{\partial \xi_i} d\xi_i \right) \end{equation*}
\end{lem}

The proof of both lemmas follows from a simple computation, just using the commutation relations mentioned above. Another important property of the multiplier $F$, easy to verify, is the following

\begin{lem} For every $a \in \S_H(\R^n) \otimes \Lambda^{\bullet} T^*\R^n$, we have 
\[ \barint [F,a] = 0 \]
\end{lem}

\subsection{Construction by excision}

The previous lemma shows that it may be relevant to consider the following cyclic cocycles on $\Psi^{-\infty}(\R^n)$, inspired of Connes' cyclic cocycles associated to Fredholm modules (see \cite{ConIHES} or \cite{ConBook}).  
\begin{equation} \label{trace cocycles}
\phi_{2k}(a_0, ..., a_{2k}) = \frac{k!}{ \i^k (2k)!} \Tr \left( a_0 [F,a_1] \ldots [F,a_{2k}] \otimes \frac{\omega^{n-k}}{n!} \right)  
\end{equation}
for $0 \leq k \leq n $. Therefore, we obtain the following result, very similar to that of Connes.

\begin{prop} The periodic cyclic cohomology classes of the cyclic cocycles $\phi_{2k}$ are independant of $k$.
\end{prop}

\begin{pr} Set 
\begin{equation} \label{transgression gamma}
\gamma_{2k+1}(a_0, \ldots , a_{2k+1}) = \frac{(k+1)!}{\i^{k+1}  (2k+2)!} \Tr \left( a_0 F[F,a_1] \ldots [F,a_{2k+1}] \otimes \frac{\omega^{n-k}}{n!} \right) 
\end{equation}
It is then a straightforward calculation to verify that $(B+b)\gamma_{2k+1} = \phi_{2k} -\phi_{2k+2}$, which shows the result. $\hfill{\square}$
\end{pr} 

At this stage, we are not very far from being done. To obtain the desired cyclic cocycles on the algebra $\S_H^0{\R^n} \otimes \Lambda^\bullet T^*\R^n $ from those previously constructed, it suffices to push the latter using excision in periodic cyclic cohomology. Indeed, as we have the pseudodifferential extension 
\begin{equation*} 
0 \to \Psi^{-\infty}(\R^n) \to \Psi^{0}_H(\R^n) \to \S_H^0{\R^n} \to 0 
\end{equation*}
we look at the image of the $(B,b)$-cocycles $\phi_{2k}$ under the boundary map 
\begin{equation*} 
\partial : \HP^0(\Psi^{-\infty}(\R^n)) \longrightarrow \HP^1(\S_H^0{\R^n}) 
\end{equation*}
Thanks to this, the cocycles (\ref{trace cocycles}) involving the operator trace, which are highly non local, will be avoided and transferred to cocycles involving the Wodzicki residue. \\

To compute the image of the the cocycles (\ref{trace cocycles}) under the excision map $\partial$, we lift the cocycles $\phi_{2k}$ on $\Psi^{-\infty}(\R^n)$ to cyclic cochains $\tilde{\phi}_{2k} \in \CC^\bullet(\Psi^{0}_H(\R^n))$ using a zeta function regularization, 
\begin{multline*}
\tilde{\phi}_{2k}(a_0, ..., a_{2k}) \\ 
= \frac{k!}{\i^k (2k)!} \frac{1}{2k+1} \sum_{i=0}^{2k} \Pf_{z=0} \Tr \left( a_0 [F,a_1] \ldots [F,a_i] \Delta^{-z/4} [F,a_{i+1}] \ldots [F,a_{2k}] \otimes \frac{\omega^{n-k}}{n!} \right)
\end{multline*}

For $k=0$, we already know that $\partial [\phi_0]$ is represented by the Radul cocycle
\begin{equation*} 
c(a_0,a_1) = \barint a_0 \delta a_1 
\end{equation*}
where $\delta a_1 = [\log \vert \xi \vert' , a_1]$. \\

Now, let $k \in \N$. Then, the usual construction of the boundary map in cohomology associated to an extension gives that $\partial [\phi_{2k}]$ is represented by the inhomogeneous $(B,b)$-cocycle 
\[  (B+b)\tilde{\phi}_{2k} = \psi_{2k-1} + \phi_{2k+1}  \in \CC^{2k-1}(\Psi^{0}_H(\R^n)) \oplus \CC^{2k+1}(\Psi^{0}_H(\R^n)) \] 
where $\psi_{2k-1} = B\tilde{\phi}_{2k}$ and $\phi_{2k+1} = b\tilde{\phi}_{2k}$ are given by 
\begin{multline} \label{psi cocycle}
\psi_{2k-1}(a_0, \ldots , a_{2k-1}) \\ 
= \frac{k!}{\i^k (2k)!} \sum_{i=0}^{2k-1} (-1)^{i+1} \barint \left( a_0 [F,a_1] \ldots [F,a_i] \delta F [F,a_{i+1}]\ldots [F,a_{2k-1}] \otimes \frac{\omega^{n-k}}{n!} \right) 
\end{multline}
\begin{multline} \label{phi cocycle}
\phi_{2k+1}(a_0, \ldots , a_{2k+1}) \\
= \frac{k!}{\i^k(2k+1)!} \sum_{i=1}^{2k+1}(-1)^{i-1} \barint \left( a_0 [F,a_1] \ldots [F,a_{i-1}] \delta a_i [F,a_{i+1}]\ldots [F,a_{2k+1}] \otimes \frac{\omega^{n-k}}{n!} \right)
\end{multline}
where we define $\psi_{-1}$ as zero. $\phi_{1}$ is precisely the Radul cocycle. For the clarity of the exposition, the calculations will be detailed later in Appendix A. Then, we have :

\begin{thm} \label{cocycles} The Radul cocycle $c$ is cohomologous in the $(B,b)$-complex, to the $(B,b)$-cocycles $(\psi_{2k-1},\phi_{2k+1})$, for all $1 \leq k \leq n$. 
\end{thm}
Indeed, usual properties of boundary maps in cohomology automatically ensures this result. As a matter of fact, one can be more precise and give explicitly the transgression cochains allowing to pass from one cocycle to another. For this, we lift the transgression cochain $\gamma$ given in (\ref{transgression gamma}) to the $(B,b)$-cochain $\tilde{\gamma} \in \CC^\bullet(\Psi_H(\R^n))$, using the same trick as before : 
\begin{multline*}
\tilde{\gamma}_{2k+1} = \dfrac{(k+1)!}{\i^{k+1} (2k+2)!} \dfrac{1}{2k+3} \left[\Pf_{z=0}  \Tr \left(a_0 \Delta^{-z/4}F  [F,a_1] \ldots  [F,a_{2k+1}] \otimes \dfrac{\omega^{n-k-1}}{n!}\right)\right. \\ 
 + \left. \left. \sum_{i=0}^{2k+1} \Pf_{z=0}  \Tr(a_0 F  [F,a_1] \ldots [F,a_{i}] \Delta^{-z/4}  [F,a_{i+1}] \ldots  [F,a_{2k+1}]  \otimes \dfrac{\omega^{n-k-1}}{n!} \right) \right]
\end{multline*}
and the term $i=0$ of the sum means $\Pf_{z=0} \Tr(a_0 F \Delta^{-z} [F,a_1] \ldots , [F,a_{2k+1}] \otimes \frac{\omega^{n-k-1}}{n!})$ in the right hand-side. 

\begin{prop} \label{transgression formulas} The inhomogeneous $(B,b)$-cochains 
\[ \tilde{\phi}_{2k} - \tilde{\phi}_{2k+2} - (B+b)\tilde{\gamma}_{2k+1} = \gamma_{2k} - \gamma'_{2k + 2} \in \CC^{2k}(\Psi^{0}_H(\R^n)) \oplus \CC^{2k+2}(\Psi^{0}_H(\R^n))  \] 
for $0 \leq k \leq n$, viewed as cochains on $\S_H(\R^n)$, are transgression cochains between $(\psi_{2k-1}, \phi_{2k+1})$ and $(\psi_{2k+1}, \phi_{2k+3})$, that is,  
\[ (\psi_{2k-1} + \phi_{2k+1}) - (\psi_{2k+1} + \phi_{2k+3}) = (B+b) (\gamma_{2k} - \gamma'_{2k + 2}) \]
Moreover, one has 
\begin{multline} \label{transgression gamma 1}
\gamma_{2k}(a_0, \ldots , a_{2k}) \\ 
= \dfrac{k!}{2 \i^{k+1} (2k+1)!} \sum_{i=0}^{2k} (-1)^i \barint \left(a_0 F  [F,a_1] \ldots [F,a_{i}] \delta F  [F,a_{i+1}] \ldots  [F,a_{2k+1}] \otimes \dfrac{\omega^{n-k-1}}{n!} \right)
\end{multline}
\begin{multline} \label{transgression gamma 2}
\gamma'_{2k}(a_0, \ldots , a_{2k}) = \barint \left(a_0 \delta a_1 [F,a_2] \ldots [F, a_{2k}] \otimes \frac{\omega^{n-k}}{n!} \right) \\
+ \dfrac{k!}{\i^{k} (2k+1)!} \sum_{i=1}^{2k} (-1)^{i-1} \barint \left(a_0 F [F,a_1] \ldots [F, a_{i-1}] \delta a_i [F, a_{i+1}] \ldots  [F,a_{2k}] \otimes \dfrac{\omega^{n-k}}{n!} \right)  
\end{multline}
\end{prop} 
That $\tilde{\phi}_{2k} - \tilde{\phi}_{2k+2} - (B+b)\tilde{\gamma}_{2k+1}$ gives a transgression cochain comes once again from the construction of a boundary map in cohomology associated to a short exact sequence. Once more, the calculations leading to these formulas are given in Appendix A. \\


\subsection{Construction with Quillen's Algebra Cochains}

The interest about Quillen's theory of cochains here is that the $(B,b)$-cocycles we want to get are obtained purely algebraically, since we do not need to pass first through $(B,b)$-cocycles on the algebra of regularizing operators. For the convenience of the reader, we briefly recall this formalism, and let him report to the original paper \cite{Qui1988} or the Appendix B for more details. 

\subsubsection{Preliminaries.} \\

Let $A$ an associative algebra over $\C$ with unit. The \emph{bar construction} $B$ of $A$ is the differential graded coalgebra $B = \bigoplus_{n \geq 0} B_n$, with $B_n = A^{\otimes n}$ for $n \geq 0$ with coproduct $\Delta : B \to B \otimes B$
\begin{equation*}\Delta (a_1, \ldots, a_n) = \sum_{i=0}^n (a_1, \ldots, a_i) \otimes (a_{i+1}, \ldots, a_n) \end{equation*}
The counit map $\eta$ is the projection onto $A^{\otimes 0} = \C$, and the differential is $b'$ :
\begin{equation*} b'(a_1, \ldots, a_{n+1}) = \sum_{i=1}^n (-1)^{i-1}(a_1, \ldots, a_i a_{i+1}, \ldots, a_{n+1}) \end{equation*}
which is defined as the zero map on $B_0$ and $B_1$. These operations confer a structure of differential graded coalgebra to $B$. \\

A \emph{bar cochain of degree n} on $A$ is a $n$-linear map over $A$ with values in an algebra $L$. These cochains form a complex denoted $\Hom(B,L)$, whose differential is given by 
\begin{equation*} \delta_{\mathrm{bar}}f = (-1)^{n+1} fb' \end{equation*}
for $f \in \Hom^n(B,L)$. Moreover, one has a product on $\Hom(B,L)$ : If f and g are respectively cochains of degrees $p$ and $q$,  it is given by
\begin{equation*} 
fg(a_1, \ldots, a_{p+q}) = (-1)^{pq} f(a_1, \ldots, a_p)g(a_{p+1}, \ldots, a_{p+q}) 
\end{equation*}
Therefore, $\Hom(B,L)$ has a structure of differential graded algebra.  \\

We next define  $\Omega^B$ and $\Omega^{B, \natural}$ to be the following bicomodules over $B$ :
\begin{align*}
& \Omega^B = B \otimes A \otimes B, \quad \Omega^{B, \natural} = A \otimes B
\end{align*}
Here, the $\natural$ in exponent means that $\Omega^{B, \natural}$ is the cocommutator subspace of $\Omega^B$. Thanks to this, one can show that the differential $\delta_\mathrm{bar}$ induced on $\Omega^{B, \natural}$ is in fact the Hochschild boundary, and deduce that the complex $(\Hom(\Omega^{B, \natural}, \C), b)$ is isomorphic to the Hochschild complex $(\CC^{\bullet}(A), b)$ of $A$, with degrees shifted by one. \\ 

We recall Quillen's terminology. Let $L$ be a differential graded algebra.  Elements of $\Hom(\Omega^B,L)$ will be called \emph{$\Omega$-cochains}, and those in $\Hom(\Omega^{B,\natural},L)$ as \emph{Hochschild cochains}. Recall also that the \emph{bar cochains} are the elements of $\Hom(B,L)$. \\

\setlength{\parindent}{6mm}
\textsc{Important fact.} 
\setlength{\parindent}{0mm} A cochain $f$ of this kind has three degrees : a $A$-degree as a multilinear map over $A$, a $L$ degree and a total degree $f$, which is sum. This is the one which will be considered. \\

The map $\natural : \Omega^{B,\natural} \to \Omega^{B}$, defined by the formula 
\begin{equation*}
\natural (a_1 \otimes (a_2, \ldots, a_n)) = \sum_{i=1}^n (-1)^{i(n-1)} (a_{i+1}, \ldots, a_n) \otimes a_1 \otimes  (a_2, \ldots, a_i) 
\end{equation*}
induces a map from Hochschild cochains to bar cochains. If we have a (graded) trace $\tau : L \longrightarrow \C$, we then obtain a \emph{morphism of complexes}
\begin{align*}
\begin{array}{c c c c}
\tau^\natural : & \Hom(\Omega^{B},L) & \longrightarrow & \Hom(\Omega^{B, \natural},\C) \\
                &  f                 & \longmapsto     & \tau^\natural(f) = \tau f \natural
\end{array}
\end{align*}

\subsubsection{Return to the initial problem.} 

We can now return to our context. Let $A$ be the algebra $\S_H^0(\R^n)$ of Heisenberg formal symbols on $\R^n = \R^p \times \R^q$, and $B$ the bar construction of $A$. Also, let $L$ be the graded algebra $\S_H^0(\R^n) \otimes \Lambda^\bullet T^*\R^n$. The product on these algebras is the star-product of symbols, twisted with the product on the exterior algebra. The injection 
\begin{equation*}
\rho : A  \longrightarrow L
\end{equation*}
is a homomorphism of algebras. As a consequence, $\rho$ should be viewed as a 1-cochain of "curvature" zero, e.g $\delta_{\mathrm{bar}} \rho + \rho^2 = 0$. We introduce a formal parameter $\epsilon$ of odd degree such that $\epsilon^2 = 0$, and shall actually work in the extended algebra 
\begin{equation*} \Hom(B,L)[\epsilon] = \Hom(B,L) + \epsilon \Hom(B,L) \end{equation*}
The role of that $\epsilon$ is to kill the powers of $\log \vert \xi \vert'$ which are not classical symbols, and to keep only its commutator with other symbols.  \\

Now, denote $\nabla = F + \epsilon \log \vert \xi \vert'$, and $\nabla^2 = F^2 + \epsilon[\log \vert \xi \vert', F]$ the square of $\nabla$, and introduce the "connection" $\nabla + \delta_{\mathrm{bar}} + \rho $. The fact that this operator does not belong to the algebra above is not a problem, since we shall only have interest in its "curvature", which is well defined, 
\begin{equation*}
K =  \nabla^2 + [\nabla, \rho] = F^2 + \epsilon[\log \vert \xi \vert' , F] + [F + \epsilon \log \vert \xi \vert' , \rho]
\end{equation*}
and its action on $\Hom(B,L)[\epsilon]$ with commutators. Here, we emphasize that the commutators involved are in fact \emph{graded commutators}. Let $\tau$ be the graded trace on $\Hom(B,L)[\epsilon] \otimes \Lambda^\bullet T^* \R^n$ given by
\begin{equation*} 
\tau(x + \epsilon y) = \barint y 
\end{equation*}

It turns out that the cocycles (\ref{psi cocycle}) and (\ref{phi cocycle}) constructed using excision in the previous section are obtained by considering the even cochain 
\begin{equation*} 
\theta = \tau(\partial \rho \cdot e^K) \in \Hom(\Omega^{B, \natural}, \C) 
\end{equation*}
where $\partial f \cdot g$ is defined, for $f,g \in \Hom(\Omega^B, L)$ of respective degrees $1$ and $n-1$, by the following formula : 
\begin{equation*} 
(\partial f \cdot g) \natural (a_1 \otimes (a_2, \ldots, a_n)) = (-1)^{\vert g \vert} f(a_1) g(a_2, \ldots , a_n)
\end{equation*}

The calculation of $\theta$ becomes easier if one remarks that 
\begin{equation*} e^K = e^{F^2} \cdot e^{[F,\rho] + \epsilon[\log \vert \xi \vert' , F + \rho]} \end{equation*} 
as $F^2 = \i \omega$ is central in $L$. Then, this easily provides that $\theta = \sum_k (\theta_{2k}' + \theta_{2k}''$), where
\begin{equation} \label{theta cochain '}
\theta_{2k}' =  \dfrac{\i^{n-k+1}}{(2k-1)!} \sum_{i=1}^{2k-1} \barint \left( \partial \rho \cdot [F,\rho]^{i-1}\delta \rho [F,\rho]^{2k-1-i} \otimes \frac{\omega^{n-k+1}}{(n-k+1)!} \right) 
\end{equation}
\begin{equation} \label{theta cochain ''}
\theta_{2k}'' = \dfrac{\i^{n-k}}{(2k)!} \sum_{i=0}^{2k-1} \barint \left( \partial \rho \cdot [F,\rho]^i \delta F [F,\rho]^{2k-1-i} \otimes \frac{\omega^{n-k}}{(n-k)!} \right)
\end{equation}

Evaluating on elements of $A$, this gives :
\begin{multline} \label{theta cochain 'ev}
\theta_{2k}'(a_0, \ldots , a_{2k-1}) \\  
= \dfrac{\i^{n-k+1}}{(2k-1)!} \sum_{i=1}^{2k-1}(-1)^{i} \barint  \left( a_0 [F,a_1] \ldots [F,a_{i-1}] \delta a_i [F,a_{i+1}]\ldots [F,a_{2k-1}] \otimes \frac{\omega^{n-k+1}}{(n-k+1)!} \right) 
\end{multline}
\begin{multline} \label{theta cochain ''ev}
\theta_{2k}''(a_0, \ldots , a_{2k-1}) \\
= \dfrac{\i^{n-k}}{(2k)!} \sum_{i=0}^{2k-1}(-1)^{i+1}\barint \left( a_0 [F,a_1] \ldots [F,a_i] \delta F [F,a_{i+1}]\ldots [F,a_{2k-1}] \otimes \frac{\omega^{n-k}}{(n-k)!} \right)
\end{multline}

The signs above not appearing in the cochains (\ref{theta cochain '}) and (\ref{theta cochain ''}) occur since the $a_i$, $\delta \rho$ and $\delta F$ are odd. \\

As announced earlier, we observe that $\theta_{2k}'$ and $\theta_{2k}''$ are up to a certain constant term the cochains $\phi_{2k-1}$ and $\psi_{2k-1}$ obtained in (\ref{psi cocycle}) and (\ref{phi cocycle}). The difference in signs is due to Quillen's formalism, which considers the total differential $B-b$, see Remark \ref{signs}.  Unfortunately, each component of $\theta_{2k} = \theta_{2k}' + \theta_{2k}''$ of $\theta$ is not a $(B,b)$-cocycle, but taking the entire cochain $\theta$ into account, this is. \\

To prove this, it only suffices to check that all the things we defined have the good algebraic properties to fit into Quillen' proof. This is the content of the following lemma, which is actually a "Bianchi identity" with respect to the "connection" $\nabla + \delta_\mathrm{bar} + \rho$.

\begin{lem} \label{bianchi identity} \emph{(Bianchi identity.)} We have $(\delta_\mathrm{bar} + \ad \rho + \ad \nabla) K = (\delta_\mathrm{bar} + \ad \rho + \ad \nabla) e^K = 0$, where $\ad$ denotes the (graded) adjoint action. \\
\end{lem}

\begin{rk} \label{rk quillen proof} The thing which guarantees this identity is that $[\nabla, \nabla] = 0$. Then, the proof is the same as that given in the paper of Quillen, \cite{Qui1988}, Section 7. Thanks to this lemma, we directly know that $(B-b) \theta = 0$, by adapting the arguments of \cite{Qui1988}, Sections 7 and 8. For the convenience of the reader, we recalled these arguments in Appendix B. This result can be refined, and we get the same results as those obtained using excision. 
\end{rk}
\begin{thm} \label{inhomogeneous theta cochains} The inhomogeneous Hochschild cochains 
\[ \theta''_{2k} - \theta'_{2k+2} \in \Hom^{2k}(\Omega^{B, \natural}, \C) \oplus \Hom^{2k+2}(\Omega^{B, \natural}, \C) \]
for $0 \leq k \leq n$, define a $(B,b)$-cocycle.
\end{thm} 

\begin{pr} Introduce a parameter $t \in \R$, and consider the following family of curvatures $(K_t)$ :
\[ K_t = \nabla^{2,t} + [tF + \epsilon \log \vert \xi \vert', \rho] \]
where $\nabla^{2,t} = F^2 + \epsilon [\log \vert \xi \vert', tF]$. Because the identity $[\nabla, \nabla^{2,t}]$ still holds, we have a Bianchi identity 
\[ (\delta_\mathrm{bar} + \ad \rho + \ad \nabla)K_t = 0 \]
Thus, the Hochschild cochain
\[ \theta^t = \tau^\natural(\partial \rho \cdot e^{K_t}) \in \Hom(\Omega^{B, \natural}, \C)[t] \] 
satisfies the relation $(B-b) \theta^t=0$ for every $t \in \R$, where we denote by $R[t]$ the polynomials with coefficients in an algebra $R$. Therefore, this relation also holds for every $k$, for the coefficient of $t^k$. This coefficient is the cochain $\theta''_{2k} + \theta'_{2k+2}$, thus, $\theta''_{2k} - \theta'_{2k+2}$ defines a $(B,b)$-cocycle.  $\hfill{\square}$
\end{pr}

Denote by $\Omega = [F, \rho] + \epsilon [\log \vert \xi \vert', \rho + F]$. The cochains which cobounds these cocycles (up to modify each of them by a constant term depending on their degrees) may be obtained rather easily by using suitable linear combinations of pairs of bar cochains $(\mu_{2j}, \mu_{2j+1})$, where $\mu$ is given by :
\begin{equation*}
\mu_{k} = \tau \left(\partial\rho \cdot \frac{e^{F^2}}{k!}\sum_{i=0}^{k} \Omega^i F \Omega^{k-i} \right)
\end{equation*}
Doing this gives transgression formulas in the spirit of those obtained in Proposition \ref{transgression formulas}.  
 
\subsection{Index theorem}

Now we know that the Radul cocycle on $\S_H^0{\R^n}$ 
\[ c(a_0, a_1) = \barint a_0 \delta a_1 \]
with $\delta a_1 = [\log \vert \xi \vert' , a_1]$, is cohomologous to the inhomogeneous $(B,b)$-cocycle 
\[ \psi_{2n-1} + \phi_{2n+1}  \in \CC^{2n-1}(\S^{0}_H(\R^n)) \oplus \CC^{2n+1}(\S^{0}_H(\R^n)) \] 
recalling that,
\begin{equation*}
\psi_{2n-1}(a_0, \ldots , a_{2n-1}) = \frac{1}{\i^n (2n)!} \sum_{i=0}^{2n-1} (-1)^{i+1} \barint  a_0 [F,a_1] \ldots [F,a_i] \delta F [F,a_{i+1}]\ldots [F,a_{2n-1}] 
\end{equation*}
\begin{multline*}
\phi_{2n+1}(a_0, \ldots , a_{2n+1}) = \\	
\frac{1}{\i^n(2n+1)!} \sum_{i=1}^{2n+1}(-1)^{i-1} \barint a_0 [F,a_1] \ldots [F,a_{i-1}] \delta a_i [F,a_{i+1}]\ldots [F,a_{2n+1}] 
\end{multline*}
it suffices to compute $\psi_{2n-1} + \phi_{2n+1}$ to obtain an index theorem. To begin, we first notice that by Lemma \ref{F differential}, we may rewrite the cocycles above as 
\begin{equation} \label{psi cocycle 1}
\psi_{2n-1}(a_0, \ldots , a_{2n-1}) =  \frac{\i^{2n-1}}{\i^{n}(2n)!} \sum_{i=0}^{2n-1} (-1)^{i+1} \barint  a_0 da_1 \ldots da_i \delta F da_{i+1} \ldots da_{2n-1} 
\end{equation}
\begin{equation} \label{phi cocycle 1}
\phi_{2n+1}(a_0, \ldots , a_{2n+1}) = \frac{\i^{2n-1}}{\i^{n}(2n+1)!} \sum_{i=1}^{2n+1}(-1)^{i-1} \barint  a_0 da_1 \ldots da_{i-1} \delta a_i da_{i+1} \ldots da_{2n+1} 
\end{equation}
The construction of the Wodzicki residue to $\Lambda^\bullet T^* \R^n$-valued symbols in the Paragraph \ref{general context} imposes that the $\littlebarint$ selects only the coefficient associated to the volume form $dx_1 \ldots dx_n d\xi_1 \ldots d\xi_n$. In (\ref{phi cocycle 1}), this coefficient must be a sum of terms of the form $\frac{\partial b_1}{\partial x_1} \ldots \frac{\partial b_n}{\partial x_n} \frac{\partial b_{2n+1}}{\partial \xi_1} \ldots \frac{\partial b_{2n}}{\partial \xi_n}$ for some Heisenberg symbols $b_1, \ldots, b_n$ of \emph{order 0}. Such terms have Heisenberg pseudodifferential order $-(p+2q)$.  \\

However, in (\ref{phi cocycle 1}), there is in each sum an additional factor of the form $\delta a_i$, which is a symbol of degree $-1$. Hence, the symbols appearing in the formula are at most of Heisenberg order $-(p+2q+1)$, and vanishes because of (\ref{wodzicki residue symbol}). \\

The formula for the cocycle (\ref{psi cocycle 1}) also reduces to a more simple one, but which is in general non-zero. A simple computation gives that 
\[ \delta F = \i \left( \sum_{i=1}^p \dfrac{\xi_i^3 d\xi_i }{\vert \xi \vert'^4} + \dfrac{1}{2}\sum_{i=p+1}^n \dfrac{\xi_i d\xi_i}{\vert \xi \vert'^4}  \right) \]   
Then, we proceed as we did to obtain the formula (\ref{phi cocycle 1}). The coefficient on $dx_1 \ldots dx_n d\xi_1 \ldots d\xi_n$ of the symbols in (\ref{psi cocycle 1}) must be of the form 
\begin{enumerate}[label={(\roman*)}]
\item $\frac{\partial b_1}{\partial x_1} \ldots \frac{\partial b_n}{\partial x_n} \frac{\partial b_{2n+1}}{\partial \xi_1} \ldots \frac{\xi_i^3}{\vert \xi \vert'^4} \ldots \frac{\partial b_{2n}}{\partial \xi_n}$ if $1 \leq i \leq p$, \\
\item $\frac{\partial b_1}{\partial x_1} \ldots \frac{\partial b_n}{\partial x_n} \frac{\partial b_{2n+1}}{\partial \xi_1} \ldots \frac{\xi_i}{\vert \xi \vert'^4} \ldots \frac{\partial b_{2n}}{\partial \xi_n}$ if $p+1 \leq i \leq n$
\end{enumerate}
where in each point, the term depending on $\vert \xi \vert'^4$ replaces the term $\frac{\partial b_{2n+i}}{\partial \xi_i}$. In all case, these terms are of order $-(p+2q)$. Thus, if we denote the Heisenberg principal symbol by
\[ \sigma : \S_H^0(\R^n) \to C^\infty(S_H^* \R^n) \]
the symbol of order $-(p+2q)$ of $a_0 da_1 \ldots da_i \delta F da_{i+1} \ldots da_{2n-1}$ is 
\begin{equation*}
\sigma(a_0) d\sigma(a_1) \ldots d\sigma(a_i) \delta F d\sigma (a_{i+1}) \ldots d\sigma(a_{2n-1}) = (-1)^i\delta F \sigma(a_0) d\sigma(a_1) \ldots d\sigma(a_{2n-1})
\end{equation*} 
We emphasize that the latter product is no more the star-product but the usual product of functions. \\

The vector field $L = \sum_{j=1}^p \xi_j \partial_{\xi_j} + 2 \sum_{j=p+1}^n \xi_j \partial_{\xi_j}$ on $T^* \R^n$ is the generator of the Heisenberg dilations. This implies that $\iota_L d\sigma(a_i) = d\sigma(a_i) \cdot L = 0$ since the $a_i$ are symbols of order $0$. Using (\ref{wodzicki residue symbol}), and observing that $\iota_L \delta F = \i$, we obtain  
\begin{equation*} \label{final cocycle}
\psi_{2n-1}(a_0, \ldots , a_{2n-1}) = \frac{(-1)^n}{(2\pi \i)^n (2n-1)!} \int_{S_H^*\R^n} \sigma(a_0) d\sigma(a_1) \ldots d \sigma(a_{2n-1}) 
\end{equation*}
So, we have proved the following theorem

\begin{thm} The Radul cocycle is $(B,b)$-cohomologous to the homogeneous $(B,b)$-cocycle on $\S_H(\R^n)$ defined by 
\begin{equation*} 
\psi_{2n-1}(a_0, \ldots , a_{2n-1}) = \frac{1}{(2\pi \i)^n (2n-1)!} \int_{S_H^*\R^n} \sigma(a_0) d\sigma(a_1) \ldots d \sigma(a_{2n-1}) 
\end{equation*}
\end{thm}

From this theorem and the pairing (\ref{index pairing}), given for any $[\phi] \in \HP^1(\S_H(\R^n)$ and $u \in K_1(\S_H(\R^n)$ by the formula 
\begin{equation*}
\langle [\phi], u \rangle = \sum_{k \geq 0} (-1)^k k! (\phi_{2k+1} \otimes \tr)(u, u^{-1}, \ldots, u, u^{-1}) 
\end{equation*}
one has the following index theorem for Heisenberg elliptic pseudodifferential operators of order $0$, which only depends on the principal symbol. Here, working in the framework of cyclic cohomology is convenient because we can directly pass from scalar symbols to matrices thanks to Morita equivalence. 

\begin{thm} Let $P \in M_N(\Psi^0_H(\R^n))$ a Heisenberg elliptic pseudodifferential operator of symbol $u \in GL_N(\S^0_H(\R^n))$, and $[u] \in K_1(\S_H^0(\R^n))$ its (odd) K-theory class. Then, we have a formula for the Fredholm index of $P$ :
\[ \Ind(P) = \Tr(\Ind [u]) = - \frac{(n-1)!}{(2\pi \i)^n (2n-1)!} \int_{S_H^* \R^n} \tr(\sigma(u)^{-1} d\sigma(u) (d\sigma(u)^{-1} d\sigma(u))^{n-1})) \]
\end{thm}

\section{Discussion on manifolds with conical singularities}

Studying index theory on manifolds with singularities is actually one of the motivations for studying a residue index formula adapted to cases where the zeta function exhibits multiple poles. It is indeed known for many years that zeta functions of differential operators on conic manifolds have double poles, see for example the paper of Lescure \cite{Lescure}. In the pseudodifferential case, even triple poles may occur, see \cite{GilLoy2002}. 

We shall first recall briefly what we need from the theory of conic manifolds, e.g pseudodifferential calculus, residues and results on the associated zeta function. This review part essentially follows the presentation of \cite{GilLoy2002}.

\subsection{Generalities on b-calculus and cone pseudodifferential operators}In our context, manifolds with conical singularities are just manifolds with boundary with an additional structure given by a suitable algebra of differential operators. \\

More precisely, let $M$ be a compact manifold with (connected) boundary, and $r : M \to \R_+$ be a boundary defining function, e.g a smooth function vanishing on $\partial M$ and such that its differential is non zero on every point of $\partial M$. We work in a collar neighbourhood $[0,1)_r \times \partial M_x$ of the boundary, the subscripts are the notations for local coordinates.  

\begin{df} A \emph{Fuchs type differential operator} $P$ of order m is a differential operator on $M$ which can be written in the form
\begin{equation*} P(r,x) = r^{-m} \sum_{j+\vert \alpha \vert \leq m} a_{j,\alpha}(r,x) (r \partial_r)^j \partial_x^\alpha \end{equation*} 
in the collar $[0,1)_r \times \partial M_x$. The space of such operators will be denoted $r^{-m} \Diff_b^m(M)$. 
\end{df}

$\Diff_b^m(M)$ are the $b$-differential operators of Melrose's calculus for manifolds with boundary. We now recall the associated small $b$-pseudodifferential calculus $\Psi_b(M)$.\\

Let $M_b^2$ be the $b$-stretched product of $M$, e.g the manifold with corners whose local charts are given by the usual charts on $M^2 \smallsetminus \partial M^2$, and parametrized by polar coordinates over $\partial M$ in $M^2$. More precisely, writing $M \times M$ near $r=r'=0$ as
\[ M^2 \simeq [0,1]_r \times [0,1]_{r'} \times \partial M^2 \]
this means that we parametrize the part $[0,1]_r \times [0,1]_{r'}$ in polar coordinates
\[ r = \rho \cos \theta , \quad r' = \rho \sin \theta \] 
for $\rho \in \R_+$, $\theta \in [0,\pi/2]$. The right and left boundary faces are respectively the points where $\theta = 0$ and $\theta = \pi/2$. \\

Let $\Delta_b$ the $b$-diagonal of $M^2_b$, that is, the lift of the diagonal in $M^2$. Note that $\Delta_b$ is in fact diffeomorphic to $M$, so that any local chart on $\Delta_b$ can be considered as a local chart on $M$.  \\ 
\begin{df} The \emph{algebra of $b$-pseudodifferential operators of order $ m $}, denoted $\Psi_b^m(M)$, consists of operators $P : C^\infty(M) \to C^\infty(M)$ having a Schwartz kernel $K_P$ such that 
\begin{enumerate}[label={(\roman*)}]
\item Away from $\Delta_b$, $K_P$ is a smooth kernel, vanishing to infinite order on the right and left boundary faces. \\

\item On any local chart of $M^2_b$ intersecting $\Delta_b$ of the form $U_{r,x} \times \R^n$ such that $\Delta_b \simeq U \times \{0\}$, and where $U$ is a local chart in the collar neighbourhood $[0,1)_r \times \partial M_x$ of $\partial M$, we have
\begin{equation*}
K_P(r,x,r',x') = \dfrac{1}{(2\pi)^n} \int e^{\i (\log(r/r') \cdot \tau + x \cdot \xi)} a(r,x,\tau,\xi) \, d\tau \, d\xi
\end{equation*} 
where $a(y,\nu)$, with $y = (r,x)$ and $\nu = (\tau, \xi)$, is a classical pseudodifferential symbol of order $m$, plus the condition that $a$ is smooth in the neighbourhood of $r=0$. 
\end{enumerate}
\end{df}

Remark that $\log(r/r')$ should be singular at $r = r' = 0$ if we would have considered kernels defined on $M^2$. Introducing the $b$-stretch product $M^2_b$ has the effect of blowing-up this singularity. \\

The \emph{algebra of conic pseudodifferential operators} is then the algebra $r^{-\Z} \Psi_b^{\Z}(M)$. The opposed signs in the filtrations are only to emphasize that $r^{\infty}\Psi_b^{-\infty}(M)$ is the associated ideal of regularizing operators. \\

To such an operator $A = r^{-p} P \in r^{-p} \Psi_b^{m}$ , we define on the chart $U$ the local density 
\begin{equation*}
\omega(P)(r,x) = \left(\int_{\vert \nu \vert = 1} p_{-n}(r,x,\tau,\xi) \iota_L d\tau d\xi \right) \cdot \, \dfrac{dr}{r} \, dx
\end{equation*}  
where $\nu = (\tau, \xi)$ and $L$ is the generator of the dilations. \\

It turns out (but this is not obvious) that this a priori local quantity does not depend on the choice of coordinates on $M$, and hence, define a globally defined density $\omega(P)$, smooth on $M$, that we call the \emph{Wodzicki residue density}. Unfortunately, the integral on $M$ of this density does not converge in general, as the boundary introduces a term in $1/r$ in the density. However, we can regularize this integral, thanks to the following lemma. Here, $\Omega_b$ denote the bundle of $b$-densities on $M$, that is, the trivial line bundle with local basis on the form $(dr/r)dx$. The following lemma from Gil and Loya is proved in \cite{GilLoy2002}.   
 
\begin{lem} \label{b-regularization} Let $r^{-p} u \in C^\infty(M,\Omega_b)$, and $p \in \R$. Then, the function
\begin{equation*}
z \in \C \longmapsto \int_M r^z u
\end{equation*}
is holomorphic on the half plane $\Re z > p$, and extends to a meromorphic function with only simple poles at $z=p, p-1, \ldots $. If $p \in \N$, Its residue at $z=0$ is given by
\begin{equation}
\Res_{z=0} \int_M r^z u(r,x) \, \frac{dr}{r} dx = \dfrac{1}{p!} \int_{\partial M} \partial_r^p (r^p u(r,x))_{r=0} \, dx 
\end{equation}
\end{lem}

Applying this regularization to the Wodzicki residue density is useful to many "residues traces" that we immediately study. 

\subsection*{Traces on conic pseudodifferential operators}

We first begin by defining different algebras of pseudodifferential operators, introduced by Melrose and Nistor in \cite{MelNis1996}. The main algebra that we shall consider is 
\begin{equation*}
A = r^{-\Z} \Psi^\Z(M) = \bigcup_{p \in \Z} \bigcup_{m \in \Z} r^{- p} \Psi^m(M)
\end{equation*}
which clearly contains the algebra of Fuchs type operators. The ideal of \emph{regularizing operators} is
\begin{equation*}
I = r^{\infty} \Psi^{-\infty}(M) = \bigcup_{p \in \Z} \bigcup_{m \in \Z} r^{- p} \Psi^m(M)
\end{equation*} 
and this explains why we note the two filtrations by opposite signs in $A$. Consider the following quotients
\begin{equation*}
I_{\sigma} = r^{\infty} \Psi^\Z(M)/I, \quad I_{\partial} = r^\Z \Psi^{-\infty}(M)/I
\end{equation*} 
Here, $I_{\sigma}$ should be thought as an extension of the algebra of pseudodifferential operators in the interior of $M$, whereas $I_{\partial}$ are     
regularizing operators up to the boundary. We finally define 
\begin{equation*}
A_\partial = \A / I_{\sigma}, \quad A_\sigma = A / I_{\partial}, \quad A_{\partial,\sigma} = A / (I_{\partial} + I_{\sigma})
\end{equation*}

\begin{df} Let $P \in r^{-p} \Psi^m(M)$ be a conic pseudodifferential operator, with $p,m \in \Z$. According to Lemma \ref{b-regularization}, define the functionals $\Tr_{\partial, \sigma}$, $\Tr_{\sigma}$ to be
\begin{gather}
\Tr_{\partial, \sigma}(P) = \Res_{z=0} \int_M r^z \omega(P)(r,x) \, \dfrac{dr}{r} dx = \dfrac{1}{p!} \int_{\partial M} \partial_r^p (r^p \omega(P)(r,x))_{r=0} \, dx \\
\Tr_{\sigma}(P) = \Pf_{z=0} \int_M r^z \omega(P) \, \dfrac{dr}{r} dx
\end{gather}
where $\Pf$ denotes the constant term in the Laurent expansion of a meromorphic function.  
\end{df}

\begin{rk} Using Lemma \ref{b-regularization}, one can show that $\Tr_{\partial, \sigma}(P)$ does not depend on the choice of the boundary defining function $r$. This is not the case for $\Tr_{\sigma}(P)$, but its dependence on $r$ can be explicitly determined, cf. \cite{GilLoy2002}. 
\end{rk}

The "Partie Finie" regularization of a trace does not give in general a trace, and this is indeed the same for the functional $\Tr_{\sigma}(P)$ acting on these algebras, the obstruction to that is precisely the presence of the boundary. However, by definition, $\Tr_{\sigma}(P)$ clearly defines an extension of the Wodzicki residue for pseudodifferential operators, one can expect that it is a trace on $I_{\sigma} = r^{\infty} \Psi^\Z(M)/I$. 

\begin{thm} \emph{(Melrose - Nistor, \cite{GilLoy2002, MelNis1996})} $\Tr_{\sigma}$ is, up to a multiplicative constant, the unique trace on the algebra $I_\sigma$
\end{thm}  

By Lemma \ref{b-regularization} and the definition above, the defect of $\Tr_{\sigma}$ to be a trace is precisely measured by $\Tr_{\partial, \sigma}(P)$, which can therefore be viewed as a restriction of the Wodzicki residue to the boundary $\partial M$. Then, the following proposition seems natural.  

\begin{thm} \emph{(Melrose - Nistor, \cite{GilLoy2002, MelNis1996})} $\Tr_{\partial, \sigma}$ is, up to a multiplicative constant, the unique trace on the algebras $A_\partial$, $A_\sigma$ and $A_{\partial,\sigma}$
\end{thm}

These two traces may be seen as "local" terms, since they only depend on the symbol of the pseudodifferential operator considered. The first can be seen as a trace on interior of $M$, the second is related to the boundary $\partial M$. There is one last trace to introduce, less easy to deal with because this one is not local. \\

Fix a holomorphic family $Q(z) \in r^{\alpha z} \Psi_b^{\beta z}(M)$, with $\alpha, \beta \in \R$, such that $Q$ is the identity at $z=0$. Take $P \in r^{-p} \Psi_b^m$, with $p,m \in \Z$ and let $(PQ(z))_\Delta$ be the restriction to the diagonal $\Delta$ of $M^2$ of the Schwartz kernel of $PQ(z)$. Melrose and Nistor noticed in \cite{MelNis1996} that $(PQ(z))_\Delta$ is meromorphic in $\C$, with values in $r^{\alpha z - p} C^\infty(M)$ with possible simple poles in the set 
\[ \left\{\dfrac{-n-m}{\beta}, \dfrac{-n-m+1}{\beta}, \ldots \right\}  \]

\begin{df} Let $P \in r^{-p} \Psi_b^m$ be a conic pseudodifferential operator. Then, we define 
\begin{equation*}
\Tr_\partial(P) = \dfrac{1}{p!} \int_{\partial M} \partial_r^p (r^p \Pf_{z=0} (PQ(z))_\Delta)_{r=0} \, dx 
\end{equation*}
If $p$ is not an integer, then, $\Tr_\partial(P)$ is defined to be $0$. 
\end{df}

\begin{rk} $\Tr_\partial(P)$ depend on the choice of the operator $Q$, but the dependence can be explicitly determined, see \cite{MelNis1996}.
\end{rk}

There is an interpretation of  $\Tr_\partial$ analogous to those of $\Tr_{\partial, \sigma}$ : If the order of $P$ is less than the dimension of $M$, then $\Tr_\partial(P)$ is a kind of $L^2$ of $P$ restricted to the boundary. This is precisely the content of the following result. 

\begin{thm} \emph{(Melrose - Nistor, \cite{GilLoy2002, MelNis1996})} $\Tr_\partial(P)$ is, up to a multiplicative constant, the unique trace on the algebra \[ I_{\partial} = r^\Z \Psi^{-\infty}(M)/I \]
\end{thm}

\subsection*{Heat kernel expansion and zeta function}

Now, let  $\Delta \in r^{-2} \Diff_b^{2}(M)$ be \emph{fully elliptic}, or \emph{parameter elliptic} with respect to a parameter $\alpha$. We refer to \cite{GilLoy2002} for the definition, what we need to know is just that this condition ensures the existence of the heat kernel $e^{-t\Delta}$ of $A$, and that operators of the type $P \Delta^{-z}$, with $P \in r^{-p} \Psi_b^m$, are of trace-class on $r^{\alpha - m}L^2_b(M)$ for $z$ in the half-plane $\Re z > max\{\frac{m + n}{2} , \frac{p}{2} \}$, $n = \dim \, M$. 

\begin{ex} \label{conic laplacian} As usual, we work in a collar neighbourhood of $M$. Then, the operator 
\begin{equation} \label{ex fully elliptic}  
\Delta = \frac{1}{r^2} \left( (r \partial_r)^2 - \Delta_{\partial M} + \frac{(n-2)^2}{4} + a^2 \right) 
\end{equation} 
where $a > 1$, is and $\alpha = 1$, is an example of such an operator. See \cite{GilLoy2002} for more details.
\end{ex}

Then, the traces introduced in the previous paragraph gives the coefficients of the expansion of $\Tr(P e^{-t\Delta})$. 

\begin{thm} \emph{(Gil - Loya, \cite{GilLoy2002})}  Under the conditions above, we have 
\begin{equation*} 
\Tr(P e^{-t\Delta}) \sim_{t \to 0} \sum_{k \geq 0} a_k t^{(k-p)/2} + (b_k + \beta_k \log t)t^{k} + (c_k + \gamma_k \log t + \delta_k (\log t)^2)t^{(k-m-n)/2}
\end{equation*}   
where 
\begin{align*}
&\beta_k = C_k (\Tr_\sigma + \Tr_\partial)(P \Delta^k) \\
&\gamma_k = C'_K \Tr_{\partial, \sigma} (P \Delta^{k-m-n}) \\
&\delta_k = C_k'' \Tr_{\partial, \sigma} (P \Delta^{k-m-n})
\end{align*}
$C_k$, $C'_K$, $C_k''$ are explicit (but not of interest for us). \\

In particular, the coefficient of $\log t$ is 
\[ -\dfrac{1}{2} \Tr_\sigma(P) - \dfrac{1}{2} \Tr_\partial(P) - \dfrac{1}{4} \Tr_{\partial, \sigma}(P) \]
and the coefficient of $(\log t)^2$ is
\[ - \dfrac{1}{4} \Tr_{\partial, \sigma}(P) \]
\end{thm}

Using a Mellin transform, we can write
\[ \Tr(P \Delta^{-z/2} = \dfrac{1}{\Gamma(z/2)} \int_0^{\infty} t^{z-1} \Tr(P e^{-t\Delta}) \, dt \]
and knowing, that $z \mapsto \int_1^{\infty} t^{z-1} \Tr(P e^{-t\Delta}) \, dt$ is entire, the asymptotic expansion of the previous proposition gives the following corollary on the zeta function. 

\begin{cor} \label{conic zeta} The zeta function $z \mapsto \Tr(P \Delta^{-z/2})$ is holomorphic in the half-plane $\Re z > max\{m + n , p \}$, and extends to a meromorphic function with at most triple poles, whose set is discrete. At $z=0$, there are simple and double poles only, which are respectively given by the terms of $\log t$ and $(\log t)^2$ in the heat kernel expansion of $\Tr(P e^{-t\Delta})$. 
\end{cor}
 
\subsection{Spectral triple and regularity} In this paragraph, we want to investigate if Fuchs type operators on conic manifolds can define an abstract algebra of differential operators, so that the local index formula we gave in the first section applies. \\

We start with a conic manifold. Let $M$ be a manifold with connected boundary, with boundary defining function $r$, endowed with the algebra of Fuchs type differential operators. The points (i), (ii), (iii) of Definition \ref{df ADO} are verified, if for example we take for $\Delta$ the fully-elliptic operator of order 2 given in Example \ref{ex fully elliptic}, and require that the order is given by the differential order. More generally, working locally in a collar neighbourhood $[0,1)_r \times \partial M_x$ of the boundary $\partial M$, elementary calculations shows that
\begin{equation} \label{commutator Fuchs order}
[r^p \Diff^m_b(M), r^{p'} \Diff^{m'}_b(M)] \subset r^{p+p'} \Diff^{m + m' -1}_b(M) 
\end{equation}
and as we shall see, the fact that the order in $r$ does not decrease is the problem. \\

Let us denote by $r^p C^\infty(\partial M)$ (find a better notation ...) the subalgebra of $C^\infty(M)$ of functions $f$ which have an asymptotic expansion 
\[ f(r,x) \sim r^p f_p(x) + r^{p+1} f_{p+1}(x) + \ldots \] 
in a neighbourhood of $r=0$. Here, the $\sim$ means that the rest of such an expansion is of the form $r^N f_N(r,x)$, with $f_N$ bounded in the collar $[0,1) \times \partial M$. The case $p=0$ actually corresponds to the smooth functions on the collar. \\ 

For the algebra of the spectral triple, it seems a good choice to look for a candidate among these classes of functions. But doing so, the formula of Lemma \ref{CM trick} is no more asymptotic in the sense of Definition \ref{asymptotic expansion}. Indeed, if $b(r,x) = r^p$ for $p \in \N$, the observation (\ref{commutator Fuchs order}) shows that the terms $b^{(k)}$ are in $r^{p-2k} \Diff_b^{k}(M)$, but by the properties of the zeta function given in the Corollary \ref{conic zeta}, the function
\[ z \longmapsto \Tr(b^{(k)}\Delta^{-k-z}) \]
is holomorphic for $\Re(z) + k > \max\left\{\frac{n+k}{2}, \frac{2k - p}{2} \right\}$, which is equivalent to $\Re(z) > \max\left\{\frac{n-k}{2}, -\frac{p}{2} \right\}$. Hence, if $p \geq 0$, the function above is in general not holomorphic at $0$ when $N$ goes to infinity. In other terms, the spectral triple we may construct will be not regular, and local index formulas of Connes-Moscovici, or those given at the beginning cannot be applied directly. As we have seen, the main problem is due to the fact that there are two notions of order : The differential order, which is local, and "the order in $r$", which is not, and comes form the presence of the boundary $\partial M$. \\

However, we may recover some interesting informations on $M$ from the zeta function. Note for instance that the higher residue $\littlebarint^2$ defined in Proposition \ref{higher WG trace} gives the trace $\Tr_{\partial, \sigma}$. $\littlebarint^1$ is, modulo some constant terms, the sum of the three functionals $\Tr_{\partial, \sigma}$, $\Tr_{\sigma}$, $\Tr_{\partial}$, which illustrates that it is no more a trace on the algebra of conic pseudodifferential operators. The next paragraph is a discussion on index theory. 

\subsection{A local index formula}

The formula of Theorem \ref{local index formula} cannot be applied directly since we are not in the context of regular spectral triples. However, there are always some relevant informations to get on index theory. \\

Let $M$ be a manifold with boundary, seen as a conic manifold, and consider the extension 
\[ 0 \to r^{\infty} \Psi_b^{-\infty}(M) \to r^{-\Z} \Psi_b^{\Z}(M) \to r^{-\Z} \Psi_b^{\Z}(M)/r^{\infty} \Psi_b^{-\infty}(M) \to 0 \] 
Here, by an \emph{elliptic pseudodifferential operator} $P \in r^{-\Z} \Psi_b^{\Z}(M)$, we shall mean that $P$ is invertible in the quotient $A = r^{-\Z} \Psi_b^{\Z}(M)/r^{\infty} \Psi_b^{-\infty}(M)$. Being \emph{fully elliptic} is an extra condition on the {indicial or normal operator}, which guarantees that $P$ is Fredholm between suitable spaces. We shall not enter into these details : What we want to investigate is just the pairing given in the paragraph (\ref{index pairing}). In particular, if $P$ is fully elliptic, then the pairing really calculates a Fredholm index. \\
 
Now, let $P,Q \in r^{-\Z} \Psi_b^{\Z}(M)$.  We can still follow the "Partie Finie" argument given in the proof of Theorem \ref{local index formula}, so that we still have the Radul cocycle
\begin{align*} 
c(P,Q) & = \Pf_{z=0} \Tr([P,Q] \Delta^{-z})  \\
& \Res_{z=0} \Tr \left(P \cdot \left( \dfrac{Q - \Delta^{-z} Q \Delta^{-z}}{z}\right)\Delta^{-z}\right) 
\end{align*} 
As we already said, the Connes-Moscovici's formula in Lemma \ref{CM trick} is no more asymptotic, but from an algebraic viewpoint, the (\ref{CM trick bis}) still holds. So, for any integer $N$, which will be thought large enough, we have  
\begin{equation*}
Q - \Delta^{-z} Q \Delta^{-z} = \sum_{k=1}^N Q^{(k)} \Delta^{-k} + \dfrac{1}{2\pi\i} \int \lambda^{-z} (\lambda - \Delta)^{-1} Q^{(N+1)} (\lambda - \Delta)^{-N-1} \, d\lambda 
\end{equation*} 
We now take advantage of the fact that the traces $\Tr_{\sigma}$ and $\Tr_{\partial, \sigma}$ vanishes when the differential order of the operators is less that the dimension of $M$. We then have the following result.

\begin{thm} Let $M$ be a conic manifold, i.e a manifold with boundary endowed with a conic metric, and let $r$ be a boundary defining function. Let $\Delta$ be the "conic laplacian" of Example \ref{conic laplacian}. Then, the Radul cocycle associated to the pseudodifferential extension  
\[ 0 \to r^{\infty} \Psi_b^{-\infty}(M) \to r^{-\Z} \Psi_b^{\Z}(M) \to r^{-\Z} \Psi_b^{\Z}(M)/r^{\infty} \Psi_b^{-\infty}(M) \to 0 \]
is given by the following \emph{non local} formula :
\begin{multline*}
c(a_0,a_1) = (\Tr_{\partial, \sigma} + \Tr_{\sigma})(a_0[\log \Delta, a_1]) - \frac{1}{2}\Tr_{\partial, \sigma}(a_0[\log \Delta,[\log \Delta, a_1]]) +  \\
+ \Tr_\partial \left( a_0 \sum_{k=1}^N a_1^{(k)} \Delta^{-k} \right) + \dfrac{1}{2\pi\i} \Tr\left(\int \lambda^{-z} a_0 (\lambda - \Delta)^{-1} a_1^{(N+1)} (\lambda - \Delta)^{-N-1}\right) \, d\lambda
\end{multline*} 
for $a_0, a_1 \in \Psi_b^{\Z}(M)/r^{\infty} \Psi_b^{-\infty}(M)$
\end{thm}
In the right hand-side, the first line consists in local terms only depending on the symbol of $P$, the second line gives the non local contributions. \\

If $P \in r^{-\Z} \Psi_b^{\Z}(M)$ is an elliptic operator, so that $P$ defines an element in the odd K-theory group $K_1^{\alg}(A)$, and $Q$ an inverse of $P$ modulo $A$, we then obtain a formula for the index of $P$. The second line of the formula above should be a part of the eta invariant (when it is defined). A perspective may be to investigate how to compare these different elements in order to get another definition of the eta invariant, suitable not only for Dirac operators but also for general pseudodifferential operators.   

\appendix

\section{Computations of Section 3.1}

We give here the details of the different computations allowing to derive the different formulas of Section 3. 

\subsection{Cocycles formulas} Recall that 
\begin{multline*}
\tilde{\phi}_{2k}(a_0, \ldots, a_{2k}) \\
= \frac{k!}{\i^k (2k)!} \frac{1}{2k+1} \sum_{i=0}^{2k} \Pf_{z=0} \Tr \left( a_0 [F,a_1] \ldots [F,a_i] \Delta^{-z/4} [F,a_{i+1}] \ldots [F,a_{2k}] \otimes \frac{\omega^{n-k}}{n!} \right)
\end{multline*}

\setlength{\parindent}{6mm}
\textsc{Formula (\ref{psi cocycle}).} We compute $\psi_{2k-1} = B\tilde{\phi}_{2k}$  
\setlength{\parindent}{6mm}
\begin{multline*} 
B\tilde{\phi}_{2k}(a_0, ..., a_{2k-1}) \\ 
= \frac{k!}{\i^k (2k)!} \frac{1}{2k+1} \sum_{i=0}^{2k} \Pf_{z=0} \Tr \left[ \left([F,a_0] \ldots [F,a_i] \Delta^{-z/4} [F,a_{i+1}] \ldots [F,a_{2k-1}] \right. \right. \\ 
 \begin{split} & - [F,a_{2k-1}] [F,a_0] \ldots [F,a_{i-1}] \Delta^{-z/4} [F,a_{i}] \ldots [F,a_{2k-2}] + \ldots \\
 & \left. \left. + (-1)^{2k-1} [F,a_1] \ldots [F,a_{i+1}] \Delta^{-z/4} [F,a_{i+2}] \ldots [F,a_{2k-1}] [F,a_0] \right) \otimes \frac{\omega^{n-k}}{n!} \right] \end{split}
\end{multline*}
Then, by the graded trace property, one can remark that all the terms of the sum $\sum_{i=0}^{2k} \ldots $ are similar, so, this sum equals $(2k+1)$ times the term $i=0$. 
\begin{multline*} 
B\tilde{\phi}_{2k}(a_0, ..., a_{2k-1}) \\ 
\begin{split}
& = \frac{k!}{\i^k (2k)!} \Pf_{z=0} \Tr \left[ \left([F,a_0] \ldots [F,a_{2k-1}]\Delta^{-z/4}  - [F,a_{2k-1}] [F,a_0] \ldots [F,a_{2k-2}] \Delta^{-z/4} \right. \right. \\
& \qquad \qquad \qquad \qquad \qquad \, \, \left. \left. + \ldots + (-1)^{2k-1} [F,a_1] \ldots [F,a_{2k-1}] [F,a_0] \Delta^{-z/4} \right) \otimes \frac{\omega^{n-k}}{n!} \right] \\
& = \frac{k!}{\i^k (2k)!} \sum_{i=0}^{2k-1} \Pf_{z=0} \Tr \left( [F,a_0] \ldots [F,a_{i}]\Delta^{-z/4}[F, a_{i+1}] \ldots [F, a_{2k-1}] \otimes \frac{\omega^{n-k}}{n!} \right)
\end{split}
\end{multline*}
where we used the graded trace property in the second equality. Then, writing $[F, a_0] = F a_0 - a_0 F$, using the fact that $F$ anticommutes with the $[F, a_i]$ and the graded trace property again, we obtain  
\begin{multline*} 
B\tilde{\phi}_{2k}(a_0, ..., a_{2k-1}) \\ 
= \frac{k!}{\i^k (2k)!} \sum_{i=0}^{2k-1} \Pf_{z=0} \Tr \left(a_0 [F, a_1] \ldots [F,a_i] ((-1)^{2k-i}\Delta^{-z/4} F - (-1)^{i} F \Delta^{-z/4}) [F,a_{i+1}] \right. \\
\shoveright{\left. \ldots [F,a_{2k-1}] \otimes \frac{\omega^{n-k}}{n!} \right)} \\
=  \frac{k!}{\i^k (2k)!} \sum_{i=0}^{2k-1} (-1)^{i+1} \Res_{z=0} \Tr \left(a_0 [F, a_1] \ldots [F,a_i] \frac{[F, \Delta^{-z/4}]}{z} [F,a_{i+1}] \right. \\
\left. \ldots [F,a_{2k-1}] \otimes \frac{\omega^{n-k}}{n!} \right)
\end{multline*} 
From Theorem \ref{local index formula}, or, to be more precise, the part of the proof allowing to pass from the Partie Finie to the residue, we finally obtain 
\begin{multline*} 
B\tilde{\phi}_{2k}(a_0, ..., a_{2k-1}) \\ 
\begin{split}
&= \frac{k!}{\i^k (2k)!} \sum_{i=0}^{2k-1} (-1)^{i+1} \barint \left( a_0 [F,a_1] \ldots [F,a_i] \delta F [F,a_{i+1}]\ldots [F,a_{2k-1}] \otimes \frac{\omega^{n-k}}{n!} \right) \\
&= \psi_{2k-1}(a_0, \ldots , a_{2k-1})
\end{split} 
\end{multline*}
$\hfill{\square}$

\setlength{\parindent}{6mm}
\textsc{Formula (\ref{phi cocycle}).}  
\setlength{\parindent}{0mm} We now compute $\phi_{2k+1} = b\tilde{\phi}_{2k}$. As $[F, \, . \,]$ is an derivation on $\S_H(\R^n)$, the following equality may be observed easily 
\begin{multline*} 
b\tilde{\phi}_{2k}(a_0, ..., a_{2k+1}) =  \frac{k!}{\i^k (2k+1)!} \sum_{i=0}^{2k} (-1)^i \Pf_{z=0} \Tr \left( a_0 [F,a_1] \ldots [F,a_i] [a_{i+1}, \Delta^{-z/4}] \right. \\ 
\left. [F,a_{i+2}] \ldots [F,a_{2k+1}]) \otimes \frac{\omega^{n-k}}{n!} \right)
\end{multline*} 
Again, from the proof of Theorem \ref{local index formula}, we finally have
\begin{multline*} 
b\tilde{\phi}_{2k}(a_0, ..., a_{2k+1}) \\
\begin{split}
&= \frac{k!}{\i^k(2k+1)!} \sum_{i=1}^{2k+1}(-1)^{i-1} \barint \left( a_0 [F,a_1] \ldots [F,a_{i-1}] \delta a_i [F,a_{i+1}]\ldots [F,a_{2k+1}] \otimes \frac{\omega^{n-k}}{n!} \right) \\
&= \phi_{2k+1}(a_0, ..., a_{2k+1})
\end{split}
\end{multline*}
$\hfill{\square}$

\subsection{Transgression formulas} We now give the details of the computations needed to obtain the formulas of Proposition \ref{transgression formulas}. Recall that 
\begin{multline*}
\tilde{\gamma}_{2k+1}(a_0, \ldots, a_{2k+1}) \\
 = \dfrac{(k+1)!}{\i^{k+1} (2k+2)!} \dfrac{1}{2k+3} \left[\Pf_{z=0}  \Tr \left(a_0 \Delta^{-z/4}F  [F,a_1] \ldots  [F,a_{2k+1}] \otimes \dfrac{\omega^{n-k-1}}{n!}\right)\right. \\ 
 + \left. \sum_{i=0}^{2k+1} \Pf_{z=0}  \Tr \left(a_0 F  [F,a_1] \ldots [F,a_{i}] \Delta^{-z}  [F,a_{i+1}] \ldots  [F,a_{2k+1}]  \otimes \dfrac{\omega^{n-k-1}}{n!} \right) \right]
\end{multline*}
where the term $i=0$ of the sum means $\Pf_{z=0} \Tr \left(a_0 F \Delta^{-z} [F,a_1] \ldots , [F,a_{2k+1}] \otimes \frac{\omega^{n-k-1}}{n!}\right)$. \\

\setlength{\parindent}{6mm}
\textsc{Formula (\ref{transgression gamma 1}).}  
\setlength{\parindent}{0mm} 
We compute $B\tilde{\gamma}_{2k+1}(a_0, \ldots, a_{2k})$. By the graded trace property, applying the operator $B$ to each term of $\tilde{\gamma}_{2k+1}$ yields the same contribution. As there are $(2k+3)$ terms, we have
\begin{multline*}
B\tilde{\gamma}_{2k+1}(a_0, \ldots, a_{2k}) = \left. \frac{(k+1)!}{\i^{k+1} (2k+2)!} \Pf_{z=0} \Tr \right(F[F,a_0] \ldots [F,a_{2k}] \\
+ \left. F[F,a_{2k}][F,a_0] \ldots [F,a_{2k-1}] + \ldots + F[F,a_{1}] \ldots F[F,a_{2k}][F,a_0] ) \Delta^{-z/4} \otimes \frac{\omega^{n-k-1}}{n!}\right)
\end{multline*}
Writing $\frac{(k+1)!}{(2k+2)!} = \frac{1}{2} \frac{k!}{(2k+1)!}  $, knowing that $F$ anticommutes to the $[F,a_i]$ and that $F^2 = \i \omega$ is central, developing $F[F,a_0]$ and finally using the graded trace property,  we obtain 
\begin{multline*}
B\tilde{\gamma}_{2k+1}(a_0, \ldots, a_{2k}) \\
= \frac{k!}{\i^{k+1} (2k+1)!} \cdot \frac{1}{2} \sum_{i=0}^{2k} \Pf_{z=0} \left((a_0 F^2 - F a_0 F) [F, a_1] \ldots \Delta^{-z/4} \ldots [F, a_{2k}]) \otimes \frac{\omega^{n-k-1}}{n!}\right)
\end{multline*}
Once again using that $F^2 = \i \omega$, we can write 
\begin{multline*}
\tilde{\phi}_{2k}(a_0, \ldots, a_{2k}) \\
= \frac{k!}{\i^{k+1} (2k+1)!} \sum_{i=0}^{2k} \Pf_{z=0} \Tr \left( a_0 F^2[F,a_1] \ldots [F,a_i] \Delta^{-z/4} [F,a_{i+1}] \ldots [F,a_{2k}] \otimes \frac{\omega^{n-k-1}}{n!} \right)
\end{multline*}
hence, 
\begin{multline*}
(\tilde{\phi}_{2k} - B\tilde{\gamma}_{2k+1})(a_0, \ldots, a_{2k}) \\
= \frac{k!}{\i^{k+1} (2k+1)!} \cdot \frac{1}{2} \sum_{i=0}^{2k} \Pf_{z=0} \left(( a_0 F^2 + F a_0 F) [F, a_1] \ldots \Delta^{-z/4} \ldots [F, a_{2k}] \otimes \frac{\omega^{n-k-1}}{n!}\right) \\
= \left. \frac{k!}{\i^{k+1} (2k+1)!} \cdot \frac{1}{2} \sum_{i=0}^{2k} \Pf_{z=0} \right( a_0 F [F, a_1] \ldots ((-1)^i F \Delta^{-z/4} - (-1)^{2k-i} \Delta^{-z/4} F) \\
\left. \ldots [F, a_{2k}] \otimes \frac{\omega^{n-k-1}}{n!}\right)
\end{multline*}
Finally, we obtain 
\begin{multline*}
(\tilde{\phi}_{2k} - B\tilde{\gamma}_{2k+1})(a_0, \ldots, a_{2k}) \\
\begin{split} 
& = \frac{k!}{2\i^{k+1} (2k+1)!} \sum_{i=0}^{2k} (-1)^i \barint \left( a_0 F [F, a_1] \ldots \delta F \ldots [F, a_{2k}] \otimes \frac{\omega^{n-k-1}}{n!}\right) \\
& = \gamma_{2k} (a_0, \ldots, a_{2k})
\end{split}
\end{multline*}
$\hfill{\square}$ \\

\setlength{\parindent}{6mm}
\textsc{Formula (\ref{transgression gamma 2}).}  
\setlength{\parindent}{0mm} 
We now calculate $b\tilde{\gamma}_{2k+1}$. Writing $a_1 F = - [F, a_1] + F a_1$ and using the derivation property of $[F, \, . \,]$, 
\begin{multline*}
b\tilde{\gamma}_{2k+1}(a_0, \ldots, a_{2k+2}) \\
= - \tilde{\phi}_{2k+2} (a_0, \ldots, a_{2k+2}) \\ 
\quad + \left. \frac{(k+1)!}{\i^{k+1} (2k+3)!} \right[ \Pf_{z=0} \left( a_0 [a_1, \Delta^{-z/4}] [F,a_2] \ldots [F, a_{2k+2}] \otimes \frac{\omega^{n-k-1}}{n!} \right) \\
+ \left. \sum_{i=0}^{2k+1} (-1)^{i} \Pf_{z=0}\left( a_0 F [F,a_1] \ldots [a_{i+1}, \Delta^{-z/4}] [F,a_2] \ldots [F, a_{2k+2}] \otimes \frac{\omega^{n-k-1}}{n!}\right)\right]
\end{multline*}
Finally, 
\begin{multline*}
(\tilde{\phi}_{2k+2} + b\tilde{\gamma}_{2k+1})(a_0, \ldots, a_{2k+2}) \\
\begin{split}
& =  \left. \frac{(k+1)!}{\i^{k+1} (2k+3)!} \right[ \barint  \left(a_0 \delta a_1 [F,a_2] \ldots [F, a_{2k+2}] \otimes \frac{\omega^{n-k-1}}{n!} \right) \\
& \qquad \qquad \qquad \qquad \quad + \left. \sum_{i=1}^{2k+2} (-1)^{i-1} \barint  \left(a_0 F  [F,a_1] \ldots \delta a_i \ldots [F, a_{2k+2}] \otimes \frac{\omega^{n-k-1}}{n!}\right)\right] \\
& = \gamma_{2k+2}(a_0, \ldots, a_{2k+2})
\end{split}
\end{multline*}
$\hfill{\square}$

\section{Complements on Section 3.2}

For the convenience of the reader, we recall here Quillen's picture of $(B,b)$-cocycles and how it is used to obtain Theorem \ref{inhomogeneous theta cochains} from the Bianchi identity of Lemma \ref{bianchi identity}. \\

\subsection{More on Quillen's formalism} Let $A$ be an associative algebra over $\C$, and $B$ be the bar construction of $A$. Recall that $\Omega^B$ and $\Omega^{B, \natural}$ are the following bicomodules over $B$ :
\begin{align*}
& \Omega^B = B \otimes A \otimes B, \quad \Omega^{B, \natural} = A \otimes B
\end{align*}

\begin{thm} \label{2 periodic complex} One has a complex of period $2$
\begin{equation*}\xymatrix{
\ldots \ar[r]^-{\overline{\partial}} & B \ar[r]^-{\beta}  & \Omega^{B, \natural} \ar[r]^-{\overline{\partial}} & B \ar[r]^-\beta & \ldots} \end{equation*}
with $\overline{\partial} = \partial \natural : \Omega^{B, \natural} \to B $, where $\natural : \Omega^{B, \natural} \to \Omega^B$, $\partial : \Omega^B \to B$, $\beta : B \to \Omega^{B, \natural}$ are defined by the following formulas :
\begin{flalign*}
&\natural (a_1 \otimes (a_2, \ldots, a_n)) = \sum_{i=1}^n (-1)^{i(n-1)} (a_{i+1}, \ldots, a_n) \otimes a_1 \otimes  (a_2, \ldots, a_i) \\
&\partial (a_{1}, \ldots, a_{p-1}) \otimes a_p \otimes  (a_{p+1}, \ldots, a_n) = (a_1, \ldots, a_n) \\
&\overline{\partial} (a_1 \otimes (a_2, \ldots, a_n)) = \sum_{i=1}^n (-1)^{i(n-1)} (a_{i+1}, \ldots, a_n, a_1, a_2, \ldots, a_i) \\
&\beta(a_1, \ldots, a_n) = (-1)^{n-1} a_n \otimes (a_1, \ldots, a_{n-1}) - a_1 \otimes (a_2, \ldots, a_n) 
\end{flalign*}
\end{thm}

As Quillen shows in \cite{Qui1988}, it turns out that the 2-periodic complex constructed above is exactly the Loday-Quillen cyclic bicomplex with degrees shifted by one, and is therefore equivalent to Connes $(B,b)$-bicomplex. The shift of the degrees makes that elements of the algebra $A$ become odd in the bar construction, while they are even in the cyclic bicomplex. \\ 

Now, let $L$ be a differential graded algebra. The maps $\overline{\partial}$ and $\beta$ of the periodic complex induces maps from bar cochains to Hochschild cochains (with values in $L$) and conversely by pull-back. The following formula is a key step.

\begin{lem} \label{quillen trace lemma} Let $f, g \in \Hom(B, L)$ be bar cochains. Then, we have
\begin{equation*} 
\beta(\tau^\natural (\partial f \cdot g)) = - \tau([f,g])
\end{equation*}
\end{lem}

We carry a purely computational proof, because of the way we introduced Quillen's formalism. A more elegant and conceptual proof is given in Quillen's article \cite{Qui1988}, paragraph 5.2. The proof of this lemma is based on the following formula, 
\begin{equation} \label{partial cochain 0}
(\partial f \cdot g) \natural (a_1 \otimes (a_2, \ldots, a_n)) = \sum_{n-p<i\leq n} (-1)^{i(n-1)} (f \cdot g) (a_{i+1}, \ldots a_n, a_1, \ldots , a_i)
\end{equation} 
where $f$ and $g$ be bar cochains of respective degrees $p$ and $n-p$. The case $p=1$ will be often used, so we give it :
\begin{equation} \label{partial cochain 1}
(\partial f \cdot g) \natural (a_1 \otimes (a_2, \ldots, a_n)) = (-1)^{\vert g \vert} f(a_1) g(a_2, \ldots , a_n)
\end{equation}

\begin{pr} Let $f$ and $g$ be bar cochains of respective degrees $p$ and $n-p$. By definition, $\beta(\tau^{\natural}(\partial f \cdot g)) = \tau (\partial f \cdot g) \natural \beta$, and using (\ref{partial cochain 0}), so, 
\begin{multline*} 
\beta(\tau^{\natural}(\partial f \cdot g))(a_1, \ldots , a_n)  \\ 
 = \tau (\partial f \cdot g) \natural (((-1)^{n-1} a_n \otimes (a_1, \ldots, a_{n-1}) - a_1 \otimes (a_2, \ldots, a_n) )  \\
 = \tau\left(\sum_{n-p<i\leq n} (-1)^{n-1}(-1)^{i(n-1)} (f \cdot g) (a_{i}, \ldots a_n, a_1, \ldots , a_{i-1})  \right. \\
    - \left. \sum_{n-p<i\leq n} (-1)^{i(n-1)} (f \cdot g) (a_{i+1}, \ldots a_n, a_1, \ldots , a_i) \right)
\end{multline*}
The first sum of the last equality can be rewritten 
\begin{multline*}
\sum_{n-p<i\leq n} (-1)^{n-1}(-1)^{i(n-1)} (f \cdot g) (a_{i}, \ldots a_n, a_1, \ldots , a_{i-1}) \\
= \sum_{n-p-1<i\leq n-1} (-1)^{i(n-1)} (f \cdot g) (a_{i+1}, \ldots a_n, a_1, \ldots , a_{i})
\end{multline*}
and noting that $(-1)^{n(n-1)} = 1$, we obtain
\begin{multline*} 
\beta(\tau^{\natural}(\partial f \cdot g))(a_1, \ldots , a_n) \\
= \tau ( (-1)^{(n-p)(n-1)}(f \cdot g)(a_{n-p+1}, \ldots , a_n , a_1 , \ldots , a_{n-p}) - (f \cdot g)(a_1, \ldots , a_n)) \\
= \tau ((-1)^{(n-p)(n-1)}(-1)^{p \vert g \vert}f(a_{n-p+1}, \ldots , a_n) g(a_1 , \ldots , a_{n-p}) - (f \cdot g)(a_1, \ldots , a_n)) \\
= \tau ((-1)^{(n-p)(n-1)}(-1)^{p \vert g \vert} (-1)^{(\vert f \vert + p)(\vert g \vert +n-p)} g(a_1 , \ldots , a_{n-p}) f(a_{n-p+1}, \ldots , a_n) \\ 
\shoveright{- (f \cdot g)(a_1, \ldots , a_n))} \\
= \tau ((-1)^{(n-p)(n-p-1)}(-1)^{\vert f \vert \cdot \vert g \vert} (g \cdot f)(a_1 , \ldots , a_{n-p}, a_{n-p+1}, \ldots , a_n)  \\
- (f \cdot g)(a_1, \ldots , a_n)  )
\end{multline*}
where we used the (graded) trace property of $\tau$ in the third equality. \\

As we have $(-1)^{(n-p)(n-p-1)} = 1$, this yields the result. $\hfill{\square}$
\end{pr}  

We can now give Quillen's picture of $(B,b)$-cocycles.  

\begin{thm} \label{Quillen cycles} Let $\theta \in \Hom(\Omega^{B,\natural}, \C)$ be a Hochschild cochain, and $\eta \in \Hom(B, \C)$ be the bar cochain defined by
\begin{equation*} 
\eta_k(a_1, \ldots, a_k) = \theta(1, a_1, \ldots, a_k)
\end{equation*}   
Suppose that for each $k$, we have
\begin{align*}
&\delta_\mathrm{bar} \eta_k = (-1)^{k} \beta \theta_{k+1}, \quad  \delta_\mathrm{bar} \theta_{k+1} = (-1)^k \overline{\partial} \eta_{k+2}
\end{align*}
and that in addition, $\theta_{n+1}(a_0,a_1, \ldots, a_n) = 0$ if $a_i = 1$, for $i \geq 1$. 

Then, for all $k$, $B \theta_{k+1} = b \theta_{k-1}$.  
\end{thm}

\begin{rk} \label{signs} This means that if we redefine signs correctly in $\theta$, we obtain a $(B,b)$-cocycle in our sign conventions. 
\end{rk}

\subsection{Complements on Remark \ref{rk quillen proof}}

We give here the details of Quillen's arguments. The only thing we have done towards the original paper \cite{Qui1988} is to mix the arguments of Sections 7 and 8.   

\begin{lem} \label{bianchi identity} \emph{(Bianchi identity.)} We have $(\delta_\mathrm{bar} + \ad \rho + \ad \nabla) K = (\delta_\mathrm{bar} + \ad \rho + \ad \nabla) e^K = 0$, where $\ad$ denotes the (graded) adjoint action. \\
\end{lem}

\begin{pr}  Let $D$ be the derivation $\delta_\mathrm{bar} + \ad \rho + \ad \nabla$. It suffices to check that $D(K) = 0$, the other equality will follow in virtue of the differentiation formula 
\begin{equation*}
D(e^K) = \int_0^1 e^{(1-s)K} D(K) e^{sK} ds
\end{equation*} 
We first remark that $[\nabla, \nabla^2] = 0$, using that $\epsilon$ commutes (in the graded sense) with elements of $\Hom(B,L)$ and that $\epsilon^2 =0$. Furthermore $\delta_\mathrm{bar} \nabla^2 = 0$ since $\delta_\mathrm{bar}$ vanishes on 0-cochains. Therefore, 
\begin{align*}
D(K) & = (\delta_\mathrm{bar} + \ad \rho + \ad \nabla)(\nabla^2 + [\nabla , \rho]) \\
     & = \delta_\mathrm{bar} [\nabla, \rho] + [\rho, [\nabla, \rho]] + [\rho, \nabla^2] + [\nabla, [\nabla, \rho]] \\
     & = [\nabla, \rho^2] + \rho [\nabla , \rho] - [\nabla , \rho] \rho + [\rho, \nabla^2] + [\nabla^2, \rho] \\
     & = 0 
\end{align*} 
The result is proved. $\hfill{\square}$
\end{pr}

According to Theorem \ref{Quillen cycles}, let us define the bar cochain $\eta \in \Hom(B,\C)$ :
\begin{equation*}
\eta_{2k-1}(a_1, \ldots, a_{2k-1}) = \theta_{2k}(1, a_1, \ldots, a_{2k+1})
\end{equation*}
Also remark that $\eta = \tau(e^K)$.

\begin{prop} The bar and Hochschild cochains $\eta$ and $\theta$ satisfies the relations 
\begin{align*}
&\delta_\mathrm{bar} \eta = \pm \beta \theta, \quad   \delta_\mathrm{bar} \theta = \pm \overline{\partial} \eta 
\end{align*} 
The $\pm$ means that the sign is positive in the even case and negative in the odd case. 
\end{prop}

\begin{pr} For the first formula of the proposition, we have
\begin{equation*}
\delta_\mathrm{bar} \eta  = \delta_\mathrm{bar} ( \tau(e^K)) = \tau (\delta_\mathrm{bar} e^K) = \tau (\delta_\mathrm{bar} e^K + [\nabla , e^K]) = -\tau ([\rho, e^K]) = \pm \beta(\tau^\natural (\partial \rho \cdot e^K)) 
\end{equation*} 
The second equality uses the trace property of $\tau$, the third is the Bianchi identity of the lemma above, and the last one is Lemma \ref{quillen trace lemma}. \\

For the second formula, first recall that $\delta_\mathrm{bar} \rho + \rho^2 = 0$. Then, one has : 
\begin{equation*}
\begin{array}{r l c l}
& \delta_\mathrm{bar} (\tau^\natural (\partial \rho \cdot e^K)) & = & \tau^\natural (\partial(-\rho^2) e^K - \partial \rho \cdot \delta_\mathrm{bar} e^K) \\
0 = & \tau^\natural([\rho, \partial \rho \cdot e^K]) & = &\tau^\natural ((\rho \cdot \partial \rho + \partial \rho \cdot \rho)e^K - \partial \rho \cdot [\rho, e^K] ) \\
0 = & \tau^\natural([\nabla, \partial \rho \cdot e^K]) & = &\tau^\natural(\partial[\nabla, \rho] e^K - \partial \rho \cdot [\nabla , e^K])
\end{array}
\end{equation*}
Adding these three equations, using Bianchi identity and $\delta_\mathrm{bar} \rho + \rho^2 = 0$ yields
\begin{equation*} 
\delta_\mathrm{bar} (\tau^\natural (\partial \rho \cdot e^K)) = \tau^\natural(\partial[\nabla, \rho] e^K) = \tau^\natural(\partial K \cdot e^K)
\end{equation*} 
The last equality follows from the definition of $K$. Moreover, 
\begin{equation*} 
\overline{\partial}(e^K) = \tau^\natural(\partial e^K) = \int_0^1 \tau^\natural(e^{(1-t)K} \cdot  \partial K \cdot e^{tK}) dt = \tau^\natural(\partial K \cdot e^K)
\end{equation*}
where last equality stands because of the trace property. 
This concludes the proof. $\hfill{\square}$
\end{pr} 

Hence, Theorem \ref{Quillen cycles} shows that $\theta$ gives rise to a $(B,b)$-cocycle (up to changing signs). The same arguments may be used to complete the proof of Theorem \ref{inhomogeneous theta cochains}.
 
\bibliographystyle{plain}
\bibliography{bibliographie}

\end{document}